\documentclass[12pt]{article}

\usepackage{amssymb,amsmath}

\def\zs#1{_{\lower5pt\hbox{$\scriptstyle#1$}}}
\def\zz#1{_{\lower4pt\hbox{$\scriptstyle#1$}}}

\newcommand{\ostar}{\textstyle\bigodot\!\!\!\!\!\!\star\,\,\,}

\newcommand{\aA}{{\mathbf A}}
\newcommand{\BB}{{\mathbf B}}
\newcommand{\CC}{{\mathbf C}}
\newcommand{\bg}{{\mathbf g}}
\newcommand{\bv}{{\mathbf v}}
\newcommand{\bs}{{\mathbf S}}
\newcommand{\bH}{{\mathbf H}}
\newcommand{\bM}{{\mathbf M}}
\newcommand{\bN}{{\mathbf N}}
\newcommand{\bL}{{\mathbf L}}

\newcommand{\bbR}{{\mathbf R}}
\newcommand{\bbS}{{\mathbf S}}

\newcommand{\cM}{{\cal M}}
\newcommand{\cD}{{\cal D}}
\newcommand{\cE}{{\cal E}}
\newcommand{\cA}{{\cal A}}
\newcommand{\cH}{{\cal H}}
\newcommand{\cN}{{\cal N}}
\newcommand{\cF}{{\cal F}}
\newcommand{\cL}{{\cal L}}
\newcommand{\cK}{{\cal K}}
\newcommand{\cP}{{\cal P}}
\newcommand{\cR}{{\cal R}}
\newcommand{\cT}{{\cal T}}
\newcommand{\cU}{{\cal U}}

\newcommand{\cX}{{\mathfrak X}}
\newcommand{\fg}{{\mathfrak g}}

\newcommand{\bR}{{\mathbb R}}
\newcommand{\bS}{{\mathbb S}}
\newcommand{\bC}{{\mathbb C}}
\newcommand{\bZ}{{\mathbb Z}}
\newcommand{\bT}{{\mathbb T}}

\newcommand{\picX}{\pi_{_{\scriptstyle\cX}}}

\newcommand{\om}{\omega}
\newcommand{\ve}{\varepsilon}
\newcommand{\of}{\overline{f}}
\newcommand{\ou}{\overline{u}}
\newcommand{\op}{\overline{p}}
\newcommand{\oz}{\overline{z}}
\newcommand{\ov}{\overline{v}}
\newcommand{\oB}{\overline{B}}
\newcommand{\oP}{\overline{P}}
\newcommand{\oU}{\overline{U}}
\newcommand{\oZ}{\overline{Z}}
\newcommand{\oD}{\overline{D}}
\newcommand{\opa}{\overline{\partial}}
\newcommand{\ochi}{\overline{\chi}}
\newcommand{\opsi}{\overline{\psi}}
\newcommand{\hb}{\hbar}
\newcommand{\pa}{\partial}
\newcommand{\wh}{\widehat}
\newcommand{\wt}{\widetilde}
\newcommand{\od}{\stackrel{\text{\rm def}}{=}}
\renewcommand{\Im}{\mathop{\rm Im}\nolimits}
\newcommand{\RR}{\mathop{\rm Re}\nolimits}
\newcommand{\rank}{\mathop{\rm rank}\nolimits}
\newcommand{\tr}{\mathop{\rm tr}\nolimits}
\newcommand{\ad}{\mathop{\rm ad}\nolimits}
\newcommand{\id}{\mathop{\rm id}\nolimits}
\newcommand{\su}{\mathop{\rm su}\nolimits}
\newcommand{\SU}{\mathop{\rm SU}\nolimits}
\newcommand{\SO}{\mathop{\rm SO}\nolimits}
\newcommand{\Hom}{\mathop{\rm Hom}\nolimits}
\newcommand{\const}{\mathop{\rm const}\nolimits}
\newcommand{\diag}{\mathop{\rm diag}\nolimits}
\newcommand{\Ker}{\mathop{\rm Ker}\nolimits}

\newtheorem{definition}{Definition}[section]
\newtheorem{lemma}{Lemma}[section]
\newtheorem{theorem}{Theorem}[section]
\newtheorem{example}{Example}[section]

\newtheorem{corollary}{Corollary}[section]

\begin{document}

\title{Quantum surfaces, special functions,\\ 
and the tunneling effect}
\author{Mikhail V. Karasev}
\date{Lett. Math. Phys., 2001, v.56, N3, 229--269}

\maketitle

\begin{abstract}
The notion of quantum embedding is considered for two classes
of examples:
quantum coadjoint orbits in Lie coalgebras
and quantum symplectic leaves in spaces with non-Lie 
permutation relations. 
A method for constructing irreducible representations 
of associative algebras
and the corresponding trace formulas over leaves 
with complex polarization are obtained.
The noncommutative product on the leaves 
incorporates a closed 2-form and a measure which (in general)
are  different from
the classical symplectic form and the Liouville measure. 
The quantum objects are related to some generalized special 
functions.  
The difference between classical and quantum geometrical
structures could even occur to be exponentially small 
with respect to the deformation parameter. 
That is interpreted as a tunneling effect in the quantum
geometry.
\end{abstract}

\bigskip
\bigskip

\hbox to \hsize{\hfill \it Dedicated to the memory of Professor M.~Flato}

\bigskip

\section{Introduction}

A manifold $\cM$ is called a {\it quantum manifold\/} 
if there is an associative noncommutative algebra of functions
over~$\cM$  
with unity element~$1$ and with involution $f\to\of$ given by
the complex conjugation.  
More precisely, following the pioneer works [2, 3], we assume that
there is a family of algebras $\cF(\cM)$  
parameterized by $\hb\geq0$ such that
the product $\star$ in $\cF$ is a deformation 
of the usual commutative product of functions:  
\begin{equation}
f\star g=fg+\hb c_1(f,g)+\hb^2 c_2(f,g)+\dots,\qquad 
\hb\to0,
\label{1}
\end{equation}
on a subspace in $\cF$ consisting of $\hb$-independent 
smooth functions. 
The coefficients $c_k$ in \eqref{1} are assumed 
to be bidifferential operators of order $k$. 

Of course, the operation 
$$
i\big(c_1(f,g)-c_1(g,f)\big)\od\{f,g\}
$$
determines the Poisson brackets over $\cM$, and so, 
the quantum manifold is automatically a Poisson manifold.

Let $\cX\subset \cM$ be one of symplectic leaves in $\cM$. Then the
Poisson structure generates on $\cX$ a symplectic form $\omega_0$ 
(see [32, 35]).
So, the leaf $\cX$ itself can be considered as a Poisson
manifold, and can be quantized. 
Denote by $*$ a quantum product in a space of functions over
$\cX$, and ask: how this product could be related to the product $\star$
over $\cM$? 

The first idea is to check whether the restriction operation
\begin{equation}
f\mapsto f\big|\zs{\cX},\qquad 
\cF_{\star}(\cM)\mapsto\cF_{*}(\cX)
\label{2}
\end{equation}
is a homomorphism of algebras. 
But in [14, 18, 48] it has been proven 
(in the case where $\cM=\fg^*$ is a Lie coalgebra and $\cX$ a
coadjoint orbit) 
that there is no quantum product $*$ on $\cX$ such that the
mapping \eqref{2} is a homomorphism. 

A way to avoid this difficulty was suggested in [26]: 
let us replace the classical restriction operation in \eqref{2}
by a {\it quantum restriction\/} operation 
\begin{equation}
f\mapsto f\big|\zs{\wh{\cX}}=f\big|\zs{\cX}+\hb\, e_1(f)
+\hb^2 e_2(f)+\dots 
\label{3}
\end{equation}
in order to preserve the homomorphism property on the quantum
level: 
\begin{equation}
f\big|\zs{\wh{\cX}}*g\big|\zs{\wh{\cX}}=(f\star g)\big|\zs{\wh{\cX}}.
\label{4}
\end{equation}

The differential operators $e_j$ in \eqref{3}
are quantum corrections to the classical restriction operation. 
They act not only along the
leaf $\cX$ but also in transversal directions, so the quantum
restriction $f\big|\zz{\wh{\cX}}$ ``feels'' not only values of $f$
on $\cX$, but also the germ of $f$ near $\cX$.
Actually, we can describe the quantum restriction
not by the formal $\hb$-power series \eqref{3}, but as the
action of a pseudodifferential operator composed with the
classical restriction: 
$$
f\big|\zs{\wh{\cX}}(x)=E\zs\cX\big(\stackrel{2}{\xi(x)},\,
-i\hb\stackrel{1}{\pa/\pa\xi}(x)\big)f(\xi(x)).
\eqno{(1.3\text{a})}
$$
Here the letters $\xi$ and $x$ designate points from $\cM$ and
$\cX$, the equation $\xi=\xi(x)$ determines the classical
embedding $\cX\subset\cM$, and the symbol $E\zz{\cX}(\xi,\eta)$ 
is
a certain smooth function on $T^{\raise3pt\hbox{$\scriptstyle*$}}_\cX\cM$.
The operators $e_j$ in \eqref{3} are obtained from (1.3a) 
by the Taylor expansion as $\hb\to0$.

An explicit calculation of the homomorphism \eqref{3}, \eqref{4}
was given in [26] for the case of products $\star,\,*$ generated
by (partial) complex structures. 
When such a homomorphism is fixed, 
we call $\cX$ a {\it quantum submanifold} of $\cM$.
The procedure \eqref{2}--\eqref{4}, (1.3a) can be called a 
{\it quantum embedding\/} of $\cX$ into $\cM$.

The simplest submanifolds are two-dimensional surfaces embedded
into Euclidean spaces. So, one could ask first of all 
about {\it quantum surfaces\/} 
homeomorphic to the plane, the sphere, the cylinder, the torus,
etc., embedded in the quantum sense into the quantized Euclidean
space $\cM=\bR^m$. 

Another interesting class of submanifolds is provided by coadjoint
orbits in Lie coalgebras $\cM=\fg^*$. The question about
quantum coadjoint orbits is very natural and attractive
remembering the negative results [14, 18, 48].

Since each surface or orbit admits not only a symplectic
structure but also a complex structure, 
it is natural to consider $*$-products on them
generated by the K\"ahlerian geometry. These are the
Wick--Klauder--Berezin products [5, 6, 33, 40, 44] 
(with some modifications). Namely, let us fix a
K\"ahler form $\omega$ and a reproducing measure $dm$ over 
the symplectic leaf $\cX$
which are certain $\hb$-deformations
of the classical form $\omega_0$ and the Liouville measure 
$dm^{\omega_0}=\frac{1}{n!}|\omega_0\wedge\dots\wedge\omega_0|$.
Then the associative product $\psi*\chi$ of two functions 
$\psi$ and $\chi$ on $\cX$ can be defined by the integral formula
\begin{equation}
(\psi*\chi)(x)=\frac1{(2\pi\hb)^n}\int_{\cX}
\psi^{\#}(x|y)\,\chi^{\#}(y|x)
\exp\bigg\{\frac{i}{\hb}\int_{\Sigma(x,y)}\omega^{\#}\bigg\}\,
dm(y),
\label{5}
\end{equation}
where $2n=\dim\cX$,
points $x,y$ are running over $\cX$,
the sign $\#$ denotes the holomorphic extension
from $\cX$ to the complexification $\cX^{\#}\approx\cX\times\cX$, 
and $\Sigma(x,y)$ denotes a quadrangle membrane in $\cX^{\#}$
whose boundary consists of paths along fibers of 
the projections 
$\cX\stackrel{\pi_-}{\longleftarrow}\cX^{\#}
\stackrel{\pi_+}{\longrightarrow}\cX$ connecting 
the points $y|y\leftarrow y|x\leftarrow x|x\leftarrow
x|y\leftarrow y|y$ (see [23, 24]). 

The measure $dm$ and the K\"ahler form $\omega$ are strongly
related to each other so that the unity function~$1$ 
should be the unity element of the product \eqref{5}; see
Section~2. 

The function algebra $\cF_*(\cX)$ with the product \eqref{5}
has a representation $\psi\to\wh{\psi}$
by the Wick pseudodifferential operators $\wh{\psi}$ acting
in the Hilbert space $\cL(\cX)$ of antiholomorphic sections
over~$\cX$; see Section~3.

The homomorphism $f\to\picX(f)$ determined by 
\begin{equation}
\picX (f)\od\wh{f\big|\zs{\wh{\cX}}}
\label{6}
\end{equation}
is an {\it irreducible Hermitian representation\/} 
of the original algebra $\cF_{\star}(\cM)$
in the Hilbert space $\cL(\cX)$. 
This representation corresponds to
the symplectic leaf $\cX\subset\cM$. 
The trace formula for this representation is
\begin{equation}
\tr\picX(f)=\frac1{(2\pi\hb)^n}\int_{\cX}\,f\big|\zs{\wh{\cX}}\,dm,
\qquad
\dim\picX=\frac1{(2\pi\hb)^n}\int_{\cX}dm.
\label{7}
\end{equation}

The operation of quantum restriction  onto the leaf can be
reconstructed from the irreducible representation 
using the formula for Wick symbols:
\begin{equation}
f\big|\zs{\wh{\cX}}(x)=\tr\big(\picX(f)\,\Pi(x)\big),\qquad
x\in\cX.
\label{8}
\end{equation}
Here $\Pi:\,\cX\to\Hom(\cL(\cX))$ is the {\it coherent mapping}
determined by
\begin{equation}
\Pi(x)^2=\Pi(x)=\Pi(x)^*,\quad 
\tr\Pi(x)=1,\quad 
\frac1{(2\pi\hb)^n}\int_{\cX}\Pi\,dm=I,
\label{9}
\end{equation}
\begin{equation}
\tr\big(\Pi(x)\Pi(y)\big)=\exp\bigg\{\frac i\hb
\int_{\Sigma(x,y)}\omega^{\#}\bigg\}.
\label{10}
\end{equation}

In particular, for any Casimir function $K$ 
(i.e., the center element in $\cF_\star(\cM)$), 
the operator $\picX(K)=\lambda\cdot I$ is scalar, 
and from \eqref{8} we see that the constant $\lambda$ is equal
to $\lambda=K\big|\zz{\wh{\cX}}$.  
So the eigenvalues  of the quantum Casimir elements are
calculated by (1.3a): 
\begin{equation}
\lambda=E\zs{\cX}\big(\stackrel{2}{\xi},-i\hb\stackrel{1}{\pa/\pa\xi}\big)
K(\xi)\big|_{\xi\in\cX}.
\label{11}
\end{equation}
\medskip

{}For example, let $\cM=\fg^*$ be a Lie coalgebra, 
and let $\cX\subset\fg^*$ be a coadjoint orbit.
The general formula \eqref{8} prompts us to define Wick symbols
of {\it group elements\/}:
$$
e^{i\eta/\hb}_*\,\od\,
e^{i\eta/\hb}\big|\zs{\wh{\cX}},\qquad \eta\in\fg.
$$
Here elements $\eta$ are considered as  
linear functions on~$\fg^*$, which can, of course, be
restricted to the orbit $\cX$.

The mapping 
$$
\exp(\eta)\mapsto e^{i\eta/\hb}_*, \qquad 
G\supset \exp(\fg)\to \cF_*(\cX),
$$
is a realization of the Lie group $G$ 
(the group corresponding to the Lie algebra~$\fg$)
in the Wick function algebra over $\cX$:
$$
e^{i\eta/\hb}_*\,*\, e^{i\eta'/\hb}_*=
e^{i\,\eta\circ\eta'/\hb}_*,\qquad 
e^{i\eta/\hb}_*\,*\, e^{-i\eta/\hb}_*=1,
$$
where $\eta\circ\eta'$ is the group (Campbell--Hausdorff) 
multiplication on~$\fg$. 

{}Formula (1.3a) applied to the exponential function reads
\begin{equation}
e^{i\eta/\hb}_*=
E\zs\cX(-i\hb\stackrel{1}{\pa/\pa\eta},\stackrel{2}{\eta})
e^{i\eta/\hb}|\zs{\cX}.
\label{12}
\end{equation}
Passing to operators, we obtain the formula 
for the irreducible representation $\pi\zs{\cX}$ of the Lie
group $G$ in the Hilbert space $\cL(\cX)$:
\begin{equation}
\picX\big(\exp(\eta)\big)=E\zs\cX
\big(-i\hb\stackrel{1}{\pa/\pa\eta},\,\stackrel{2}{\eta}\big)
\widehat{e^{i\eta/\hb}|\zs{\cX}}.
\label{13}
\end{equation}

At last, using the general trace formula \eqref{7}, 
we obtain the {\it character of the irreducible representation\/}: 
\begin{align}
\tr\picX\big(\exp(\eta)\big)
&=\frac1{(2\pi\hb)^{\dim\cX/2}}
\int_{\cX}e^{i\eta/\hb}_*\,dm\nonumber\\
&=\frac1{(2\pi\hb)^{\dim\cX/2}}E\zs\cX
\big(-i\hb\stackrel{1}{\pa/\pa\eta},\,\stackrel{2}{\eta}\big)
\int_{\cX}e^{i\eta/\hb}\,dm.
\label{14}
\end{align}
The measure $dm$ in this case is proportional to the classical 
Liouville measure $dm^{\omega_0}$ on the orbit. 
The last integral looks very similar to the Kirillov character
formula [31], but the operator of the quantum restriction onto
the orbit~$\cX$
(presented by the symbol $E\zz{\cX}$) and the group elements 
$e^{i\eta/\hb}_*$ are  certainly new objects in this framework. 

Note that if the orbit $\cX$ admits 
a $G$-invariant complex structure, 
then the symbol $E\zz{\cX}$ and the functions 
$e^{i\eta/\hb}_*$ are very easily and explicitly calculated
just by solving a certain first order differential equation over
$\cX$ or by evaluating the area of certain membranes in
$\cX\times\cX$; see Section~4.
In this way we obtain the explicit formula \eqref{13} 
for irreducible representations of the Lie group by Wick
pseudodifferential operators over coadjoint orbits 
and the expression \eqref{14} for characters of these
representations, as well as the eigenvalues~\eqref{11} of the
Casimir elements.
\bigskip

In Section~5 we leave the case of Lie algebras. We consider 
general function algebras with complex 
polarization and demonstrate explicit formulas for their
Hermitian irreducible representations. A notion of {\it special
functions\/} is associated with quantum K\"ahler leaves. 

Then we are concentrated on the second example: 
$2$-dimensional surfaces of revolution $\cX\subset\bR^m$ endowed
with a generic (not group-invariant) complex structure. 
The quantum realization of such surfaces is described in
Section~6 using the results of Section~5, and we focus the
attention on the algebraic polynomial case.

Assume that the surface $\cX$ is embedded into $\cM=\bR^m$ by
means of algebraic equations, and moreover, the complex
structure on $\cX$ is taken in such a way that all 
the operators $\picX(f)$  are {\it differential\/} 
(i.e., not generic pseudodifferential) for polynomial $f$.
It is always possible to realize this situation 
in the case where $\cX$ is homeomorphic to the plane or the sphere.
Then the quantum K\"ahler form $\omega$ 
and the measure $dm$ in \eqref{5} are very special; 
namely, in local complex coordinates $z$ we have
\begin{equation}
\omega=i\hb\opa\pa\ln k(|z|^2),\qquad
dm=k(|z|^2)\ell(|z|^2)d\oz dz,
\label{15}
\end{equation}
where $k$ and $\ell$ are certain {\it hypergeometric functions}.
And vise versa: any hypergeometric function is related in this
way to some quantum surface homeomorphic to the plane or sphere. 
This was proved in [29]; see Section~6.

In the case where $\cX$ is homeomorphic to the cylinder,
the representation $\picX(f)$ cannot act by purely
differential operators, but must also include the shift operator
$\exp\{\hb\pa/\pa\oz\}$.
The corresponding $*$-product is given again by \eqref{5} and
\eqref{15}, where the functions $k$ and $\ell$ are some 
{\it theta-functions}, and the argument $|z|^2$ in \eqref{15}
must be replaced by $z+\oz$.

Note that our way to relate the theta-function to the cylinder
is different from that 
used in [47] for the construction of the Weyl quantization over
the torus, 
as well as from the approaches based on discrete subgroups 
of the Heisenberg group in [15,~36].

What is important and unexpected is the $\hb$-expansion
of the quantum K\"ahler form $\omega$ and the reproducing
measure $dm$ in the cylindric case. 
We prove in Section~7 that these quantum objects differ from the
classical form $\omega_0$ and the classical Liouville measure 
$dm^{\omega_0}$ on the symplectic leaf $\cX$ 
by exponentially small quantum corrections of order 
$O(e^{-\pi^2/\hb})$.
These corrections are given precisely by some theta-series, and
they control the difference between the topology of the cylinder
and the topology of the plane. The same results hold in the case 
of torus~[30].

\setcounter{equation}{0}

\section{Reproducing measure}

{}First of all we introduce and discuss some general definitions
related to the Wick quantization of K\"ahlerian manifolds. 

Let $\cX$ be a K\"ahlerian manifold with the K\"ahler form
$\omega$. 
{}Following [16] consider a Hermitian line bundle over $\cX$
whose curvature is $i\omega/\hb$. The Hermitian bilinear form on
each fiber is $(u,v)_\hb=u\ov\exp(-F/\hb)$, where $F$ is a local
K\"ahler potential, i.e., 
\begin{equation}
\omega=i\opa\pa F.
\label{2.1}
\end{equation}
If the cohomology class of the curvature is multiple $2\pi i$,
i.e.,  
\begin{equation}
\frac1{2\pi\hb}[\omega]\in H^2(\cX,\bZ),
\label{2.2}
\end{equation}
then this line bundle admits  nontrivial antiholomorphic
sections. 
Each section $u$ has its Hermitian fiberwise norm
$\rho_u=(u,u)_\hb$ called a {\it density function}. 

Denote by $\cL=\cL(\cX)$ the space of antiholomorphic sections with 
Hilbert norm
\begin{equation}
\|u\|\zs{\cL}=\bigg(\frac1{(2\pi\hb)^n}\int_{\cX}\rho_u\,dm\bigg)^{1/2}
\label{2.3},
\end{equation}
where $2n=\dim\cX$ and $dm$ is a smooth positive measure on $\cX$.

The space $\cL$ can be characterized by its reproducing
kernel [9, 45] 
$\cK=\sum_j|u_j|^2$, \ 
$\{u_j\}$ is an orthonormal basis in $\cL$.
Of course, $\cK$ depends not only on the choice of the K\"ahler
structure over $\cX$ but also on the choice of the measure $dm$.
The reproducing kernel determines a new K\"ahler form 
$\omega_m=i\hb\opa \pa\ln\cK$.
So, we obtain a transform $\omega\to\omega_m$ of K\"ahlerian
structures over~$\cX$. 

\begin{definition}
\rm
The measure $dm$ on $\cX$ is called a {\it reproducing measure} 
corresponding to the K\"ahler form $\omega$ if
$\omega_m=\omega$. 
\end{definition}

In general, the question about the existence and uniqueness of
the reproducing measure corresponding to the given form~$\omega$
is open. 

{}For an arbitrary measure $dm$ let us define the function
\begin{equation}
\eta=\cK\exp(-F/\hb),
\label{2.4}
\end{equation}
and introduce the new form and the new measure
\begin{equation}
\omega'=\omega+i\hb\opa \pa\ln\eta,\qquad dm'=\eta dm.
\label{2.5}
\end{equation}
Under this replacement the scalar product \eqref{2.3} and so the
space $\cL$ and the reproducing kernel $\cK$ are not changed. 

\begin{lemma}
The measure $dm'$ is the reproducing measure corresponding to
the K\"ahler form $\omega'$.
\end{lemma}

Thus we see that for each K\"ahler form there is always another 
K\"ahler form (at the same cohomology class) for which the
reproducing measure does exist.

Note that the natural choice of the measure on $\cX$ is the 
Liouville measure 
$dm^\omega=\frac1{n!}|\omega\wedge\dots\wedge\omega|$. 
{}For this choice the function $\eta$ \eqref{2.4} was first
introduced in [43], and the question of whether $\eta=\const$ or
not was raised. 
If $\eta=\const$ then $\omega'=\omega$ in \eqref{2.5} and the
reproducing measure corresponding to $\omega$ is 
$dm'=\const\cdot dm^\omega$. 
{}For instance, on homogeneous K\"ahlerian manifolds this is the
case (see [8, 43]). 

In what follows we {\it assume that the reproducing measure over
$\cX$ exists}. 

\smallskip

{}For each $x\in\cX$ denote by $\Pi(x)$ the linear operator  
in $\cL$ defined by 
\begin{equation}
(\Pi(x)u,u\big)_{\cL}=\rho_u(x),\qquad 
\forall u\in\cL.
\label{2.6}
\end{equation}

\begin{lemma}
{}For the mapping $\Pi:\,\cX\to\Hom\cL$,  
all properties \eqref{9} hold. 
In particular, we have the estimate
\begin{gather*}
\rho_u(x)\leq\|u\|^2\zs{\cL},\qquad \forall u\in\cL,\quad 
\forall x\in\cX.
\end{gather*}
\end{lemma}

Using the terminology [23], we call
\begin{equation}
p(x,y)=\tr\big(\Pi(x)\Pi(y)\big)
\label{2.7}
\end{equation}
the {\it probability function}. Its properties are the following:
$$
0\leq p(x,y)\leq 1,\qquad p(x,x)=1,
$$
\begin{equation}
\frac1{(2\pi\hb)^n}\int_{\cX}p(x,y)\,dm(y)\equiv 1,\qquad
\forall\,x\in\cX. 
\label{2.8}
\end{equation}
In the exponential representation
$p(x,y)=\exp\{-d(x,y)^2/2\hb\}$
the function $d(\,\cdot,\cdot\,)$ is the Calabi distance
between points of $\cX$ (for details see [13]).

An important question about the Calabi distance: 
does it actually distinguish points of $\cX$, 
that is, 
does it follow from $d(x,y)=0$ that $x=y$?
Or, using another language: 
\begin{equation}
p(x,y)=1\quad \Longrightarrow\quad x=y.
\label{2.9}
\end{equation}
If this property holds, the manifold $\cX$ is called the
{\it probability space\/} (see in [4, 39]).
Some sufficient condition for \eqref{2.9} was mentioned in [13].
In what follows, we assume that property \eqref{2.9} holds. 

Now, following [20, 21, 46], consider the {\it
complexification\/}  $\cX^{\#}=\cX\times\cX$ of the manifold
$\cX$.  
There are projections $\pi_+$ and $\pi_-$ from $\cX^{\#}$ to the
multipliers $\cX$. Points in $\cX^{\#}$ we denote by $x|y$, 
so that $\pi_+(x|y)=y$, $\pi_-(x|y)=x$.
Also let us identify $\cX$ with the diagonal 
$\diag(\cX\times\cX)$, i.e., $x\equiv x|x$.
The holomorphic differential on $\cX^{\#}$ is defined by 
$\pa^{\#}_{x|y}=(\opa_x,\pa_y)$,
and the tangent space $T_{x|x}\cX^{\#}$ is identified with the
direct sum of polarizations 
$T^{(0,1)}_{x}\cX\oplus T^{(1,0)}_{x}\cX$, 
so that fibers if $\pi_+$ and $\pi_-$ are considered
as integral leaves of the complex polarization and its
conjugate. 

Note that the reproducing kernel $\cK$ is naturally extended to
the complexification:
$\cK^{\#}(x|y)=\sum_ju_j(x)\overline{u_j(y)}$ to be
$\pa^{\#}$-holomorphic over $\cX^{\#}$. 
Also the form 
$\omega^{\#}_{x|y}=i\hb\opa_x\pa_y\ln p(x,y)
=i\hb\opa\pa\ln\cK^{\#}(x|y)$
is the holomorphic extension of the K\'ahler form $\omega$ 
from $\cX$ to the complexification $\cX^{\#}$.
The set $\Sigma_{\#}$ of singularities of $\omega^{\#}$ consists
of all pairs of points $x,y\in\cX$ for which $p(x,y)=0$ 
(or the Calabi distance $d(x,y)=\infty$).
The dimension of $\Sigma_{\#}$ does not exceed 
$2n=\dim\cX$, and intersections of $\Sigma_{\#}$
with fibers of $\pi_+$ and $\pi_-$ are transversal. 

Thus we can integrate
the closed form $\omega^{\#}$ over two-dimensional membranes in
$\cX\times\cX$ whose boundaries belong to fibers of $\pi_+$ or
$\pi_-$. 
Note that $\omega^{\#}$ vanishes on these fibers and so the
integral of $\omega^{\#}$ does not depend on the shape of a
boundary-path along fibers of $\pi_+$ or $\pi_-$.
This integral also does not depend on the shape of the membrane 
itself because of closedness of the form $\omega^{\#}$ and
because of the quantization conditions \eqref{2.2} around
two-dimensional holes.  

Let us specify the shape of the membrane.
Take two point $x,y\in\cX$ and identify them with points $x|x$
and $y|y$ on $\diag(\cX\times\cX)$. 
Then consider points $x|y$ and $y|x$ in $\cX^{\#}$ which are
points of intersection of $\pi_+$-fibers with $\pi_-$-fibers
over $x$ and $y$. 
Take any closed path with pieces running along $\pi_+$- and
$\pi_-$-fibers  
$x|x\gets x|y \gets y|y \gets y|x \gets x|x$
and consider a quadrangle membrane $\Sigma(x,y)$ 
in $\cX^{\#}$ whose boundary coincides with this closed
path. 
Then we obtain formula \eqref{10}
\begin{equation}
p(x,y)=\exp\Big\{\frac i\hb\int_{\Sigma(x,y)}\omega^\#\Big\}.
\label{2.10}
\end{equation}

Thus the probability function is determined by the K\"ahler form
only and, in particular, it does not explicitly depend on the
reproducing measure. 
This fact allows us to consider the property \eqref{2.8}
together with \eqref{2.10} as a {\it linear equation for the
reproducing measure\/} corresponding to the K\"ahler
form.   

{}From this equation it is possible, for instance, to find the
formal $\hb$-power expansion of the reproducing measure
under the {\it a priori\/} assumption that the measure is
$\hb$-smooth as $\hb\to0$.
Namely, the integral in \eqref{2.8} can be asymptotically
evaluated by the stationary phase method as in [13, 43]. 
Under condition \eqref{2.9} there is a unique point $y=x$ where
the phase $\int_{\Sigma(x,y)}\omega^{\#}$ takes its minimum (zero)
value, and this minimum is not degenerate. 
So, if we represent the unknown measure $dm$ in \eqref{2.8}
by some function $\sigma$ as follows:
\begin{equation}
dm=\sigma\,dm^\omega,\qquad\text{where}\quad 
dm^\omega=\frac1{n!}|\omega\wedge\dots\wedge\omega|,
\label{2.11}
\end{equation}
and introduce Hermitian matrices
$(\!(\omega_{\nu\mu})\!)$ and $(\!(\omega^{-1\mu\nu})\!)$ 
of the K\"ahler form and of the Poisson tensor 
with respect to  local complex coordinates on~$\cX$:
$$
\omega=i\omega_{\nu\mu}(x)d\oz^\nu(x)\wedge dz^\mu(x),\qquad 
\{z^\mu,\oz^{\nu}\}=i\omega^{-1\mu\nu},
$$ 
then equation \eqref{2.8} becomes asymptotically equivalent to 
\begin{equation}
\frac1{(2\pi\hb)^n}\int_{\bC^n}d\ov dv 
\exp\Big\{-\frac1\hb\Omega v\ov+(\ov\opa+v\pa)\Big\}
(\sigma\det\omega)=1.
\label{2.12}
\end{equation}
Here $\pa=\pa/\pa z(x)$ act first, 
the matrix $\Omega$ is the deformation of $\omega$: 
$$
\Omega_{\nu\mu}(x,\bv)=|\varphi(v\pa)|^2\omega_{\nu\mu}(x),
$$
where $\varphi(\xi)\,\od\,(e^\xi-1)/\xi$, and $\bv$ are complex
coordinates of vectors $V\in T_x\cX\approx \bC^n$.  
The explicit derivation of the Gaussian integral in
\eqref{2.12} transforms this equation to the following one
(derivatives act first):
\begin{equation}
f_\hb(x,\nabla)\sigma=1,\qquad 
\nabla\,\od\,\frac1{\det\omega(x)}\cdot d_x\cdot\det\omega(x),
\label{2.13}
\end{equation}
where $f_\hb$ is the function on $T^*\cX$ given by
\begin{equation}
f_\hb(x,p)=\det\omega(x)
\bigg[\frac{\exp\{\hb\Omega^{-1}(x,\cdot)p\,p\}}
{\det\Omega(x,\cdot)}\bigg]^\flat.
\label{2.14}
\end{equation}
Here the contraction operation $\flat$ applied to any function
$g(V,p)$ on $T_x\cX\times T^*_x\cX$ polynomial in~$p$ 
produces a polynomial function on $T^*_x\cX$ by the formula
$$
g(\cdot,p)^\flat \,\od\, g(\stackrel{2}{d_p},
\lower1pt\hbox{$\stackrel{1}{p}$})1, 
\qquad p\in T^*_x\cX.
$$
Explicit expression for function \eqref{2.14} is 
\begin{align}
&f_\hb=\sum_{k\geq0}\frac1{k!}
\sum_{
\begin{smallmatrix}
|\alpha_0|,|\beta_0|\geq0\\
|\alpha_1|+|\beta_1|\geq1\\
\dotfill\\
|\alpha_k|+|\beta_k|\geq 1
\end{smallmatrix}
}
\frac{\hb^{|\alpha_0|+\dots+|\alpha_k|}}
{(|\alpha_1|+1)\dots(|\alpha_k|+1)(|\beta_1|+1)\dots(|\beta_k|+1)
\alpha!\beta!}\nonumber\\
&
\cdot\!\det\omega\!\cdot\!
\opa^{\alpha_1}\pa^{\beta_1}\om_{\nu_1\mu_1}\!\dots 
\opa^{\alpha_k}\pa^{\beta_k}\om_{\nu_k\mu_k}
\frac{\pa^k}{\pa\om_{\nu_1\mu_1}\!\dots \pa\om_{\nu_k\mu_k}}
\bigg(\frac{\cF^{\alpha_0+\dots+\alpha_k}_{\beta_0+\dots+\beta_k}
(\om^{-1})}{\det\om}\bigg) \op^{\alpha_0}p^{\beta_0}\!.
\label{2.15}
\end{align}
Here we use the following notations: for arbitrary matrix~$A$
and multi-indices $\alpha, \beta$ 
the polynomial $\cF^\alpha_\beta(A)$ in matrix elements of~$A$
is given by $\cF^\alpha_\beta(A)\,\od\,
(A\pa/\pa \bv)^\alpha \bv^\beta\big|_{\bv=0}$.  

As it follows from \eqref{2.14}, \eqref{2.15}:
\begin{equation}
f_\hb(x,p)=1+\hb f^{(1)}(x,p)+\hb^2 f^{(2)}(x,p)+\dots,
\label{2.16}
\end{equation}
where $f^{(k)}$ are explicitly given symbols 
on $T^*\cX$ polynomial in $p$ of degree $2k$. 
The corresponding differential operators $f^{(k)}(x,\nabla)$ are
interesting geometric invariants of the K\"ahlerian structure 
$\omega$ (see [26]). For instance,
$$
f^{(1)}(x,\nabla)=\frac12(\Delta+\omega^{-1\mu\nu}\rho_{\nu\mu}),
$$
where $\Delta$ is the Laplace--Beltrami operator on $\cX$ and  
$\rho_{\nu\mu}=\opa_\nu\pa_\mu\ln\det\omega$ is the matrix of
the Ricci form  
$$
\rho=i\rho_{\nu\mu} d\oz^\nu\wedge dz^\mu.
$$

The formal asymptotic solution of equation \eqref{2.13} is the
following 
\begin{equation}
\sigma=1+\sum_{s\geq 1}\hb^s\sigma_s,\qquad 
\sigma_s=\sum^s_{l=1}\,\,\sum_{\substack{k_1,\dots,k_l\geq1\\|k|=s}}
(-1)^lf^{(k_1)}(x,\nabla)\dots f^{(k_l)}(x,\nabla)1.
\label{2.17}
\end{equation}
In particular, 
$$
\sigma_1=-\frac12 \omega^{-1\mu\nu}\rho_{\nu\mu}
$$
is half the scalar curvature of the  K\"ahler manifold $\cX$.

\begin{theorem}
Let the K\"ahler form $\om$ over $\cX$ satisfy conditions
\eqref{2.2}, \eqref{2.9}.
Then the formal $\hb$-expansion for the reproducing measure over
$\cX$ is given by  
\begin{equation}
dm\sim(1+\hb\sigma_1+\hb^2\sigma_2+\dots)\,dm^\om,
\label{2.18} 
\end{equation}
where $dm^\om$ is the Liouville measure \eqref{2.11},  
functions $\sigma_s$ are determined by \eqref{2.17}, and
symbols $f^{(k)}$ are taken from \eqref{2.15}, \eqref{2.16}. 
\end{theorem}

Now from \eqref{9} we obtain a useful corollary.

\begin{corollary}
Under conditions \eqref{2.2}, \eqref{2.9},
the dimension of the space of antiholomorphic sections over the
compact K\"ahler manifolds $\cX$ can be calculated by the formula  
\begin{equation}
\dim\cL(\cX)=\frac1{(2\pi\hb)^n}\int_{\cX}
(1+\hb\sigma_1+\dots+\hb^n\sigma_n)\,dm^\om,
\label{2.19}
\end{equation}
where the functions $\sigma_s$ are given in \eqref{2.17}.
\end{corollary}

Note that \eqref{2.19} is a precise formula, not asymptotic
in~$\hb$, 
and it contains only the first $n$ coefficients $\sigma_s$ 
($1\leq s\leq n$).
Indeed,  in all higher terms (for $s>n$) the integrals 
$\frac{\hb^s}{(2\pi\hb)^n}\int_{\cX}\sigma_s\,dm^\om=O(\hb^{s-n})$
are asymptotically small, and so, they should be just zero,
since the left-hand side of \eqref{2.19} is integer. 

Also note that each $\sigma_s$ in \eqref{2.19} is a functional
of the K\"ahler metric $\om$ homogeneous of degree $-s$, and so,
if we include $\hb$ 
into the definition of the K\"ahler form on $\cX$, then $\hb$
disappears from \eqref{2.19}. This means that one can set $\hb=1$
simultaneously in \eqref{2.19} and in \eqref{2.2}.

An analogous way to derive formulas of type \eqref{2.19} 
(but not in this explicit form) 
was suggested in [41];
for instance, in [41] there is a discussion about relations 
between $*$-product and the Riemann--Roch--Hirzebruch theorem.
Our formulas for $\sigma_s$ on the right-hand side of
\eqref{2.19} just represent in some way the
Riemann-Roch number and Hilbert--Samuel polynomial [11, 12, 38]. 

In particular, for compact $2$-dimensional surfaces $\cX$ 
we conclude from the Gauss--Bonnet theorem that the integral 
$\frac{1}{\pi}\int_{\cX}\sigma_1\,dm^\omega$ is the Euler number
$\chi(\cX)=c_1(\cX)$
of the surface. Thus in this case formula \eqref{2.19} reads: 
$\dim\cL(\cX)=N+\frac12\chi(\cX)$, 
where $N=\frac1{2\pi\hb}\int_{\cX}|\omega|\in\bZ_+$.

\setcounter{equation}{0}

\section{Quantization by complexification}

Over any K\"ahlerian manifold $\cX$
the probability function \eqref{2.10} determines the probability
operator $\cP$ acting by the formula
\begin{equation}
(\cP\psi)(x)\od\frac1{(2\pi\hb)^n}\int_{\cX}p(x,y)\psi(y)\,dm(y).
\label{3.1}
\end{equation}

\begin{lemma}
The probability operator is a positive self-adjoint contraction
in the space $L^2(\cX,dm)$.
\end{lemma}

We denote by $M=M(\cX)$ the range of the operator $\cP$. Then
$L^2=M\oplus\Ker\cP$.
The space $M$ is $\cP$-invariant and endowed with the norm 
\begin{equation}
\|\psi\|\zs{M}=(\cP\psi,\psi)^{1/2}\zz{L^2}.
\label{3.2}
\end{equation}
Each function $\psi\in M$ can be represented as
$\psi=\cP\varphi$, where $\varphi\in M$, and so we can define
\begin{equation}
\psi^{\#}(x|x')\od\frac1{(2\pi\hb)^n}\int_{\cX}
\exp\bigg\{\frac{i}{\hb}\int_{\Sigma(x|x',y)}\omega^{\#}\bigg\}
\varphi(y)\,dm(y),
\label{3.3}
\end{equation}
where $\Sigma(x|x',y)$ is a membrane in $\cX^{\#}$ with the
boundary $y|y\gets y|x'\gets x|x'\gets x|y\gets y|y$.
The function $\psi^{\#}$ \eqref{3.3} is the
$\pa^{\#}$-holomorphic extension of $\psi$ to the
complexification $\cX^{\#}$. 
Of course, $\psi^{\#}(x|x)=\psi(x)$.

Now let us introduce the Hilbert norm
\begin{equation}
\|\psi\|\zs{W}=\bigg(\frac1{(2\pi\hb)^{2n}}\iint_{\cX\times\cX}
|\psi^{\#}(x|x')|^2 p(x,x')\,dm(x)dm(x')\bigg)^{1/2}
\label{3.4}
\end{equation}
and denote by $W=W(\cX)$ the completion of the space $M$ by
this norm. 
The probability operator is an isometry
$M\stackrel{\cP}{\longrightarrow}W$.

{}For each $\psi\in W$ one can define the 
{\it Wick pseudodifferential operator\/} $\wh{\psi}$ in the
space $\cL$ by the bilinear form 
\begin{equation}
(\wh{\psi}u,u)\zs{\cL}\od(\psi,\rho_u)\zs{W},\qquad \forall u\in\cL;
\qquad\text{or}\quad\wh{\psi}=(\psi,\Pi)\zs{W}.
\label{3.5}
\end{equation}
Explicit formula is the following:
\begin{equation}
(\wh{\psi}u)(x)=\frac1{{(2\pi\hb)}^n}\int_{\cX}
\cK^{\#}(x|y)\psi^{\#}(x|y)u(y)e^{-F(y|y)/\hbar}\,dm(y).
\label{3.6}
\end{equation}
The function $\psi$ is called the low symbol [33], 
or the Wick symbol, or the covariant symbol [6]
of the operator $\wh{\psi}$. 
The symbol is reconstructed by the formula
\begin{equation}
\psi(x)=\tr\big(\wh{\psi}\Pi(x)\big),
\label{3.7}
\end{equation}
and so, the correspondence $\psi\to\wh{\psi}$ is one-to-one. 
Of course, the complex conjugate function corresponds to the
adjoint operator: $\overline{\psi}\to\wh{\psi}^*$.
Pure states are exactly Wick operators corresponding to density
functions:
$$
\wh{\rho}_u\bv=(\bv,u)\zs{\cL}u.
$$
Positive operators have positive symbols. 

The reproducing kernel is the ``eigenfunction'' of all Wick
operators: 
$$
\wh{\psi}\cK^{\#}=\psi^{\#}\cK^{\#}
$$
(where $\wh{\psi}$ acts by the first argument, i.e., 
$\wh{\psi}\approx\wh{\psi}\otimes I$).

Now let us introduce an associative multiplication to the space
of Wick symbols. Note that there is a natural ``matrix'' or
groupoid multiplication of sections generated by the scalar
product \eqref{2.3}, namely, 
$$
(\psi\times\chi)(x)=\frac1{(2\pi\hb)^n}\int_{\cX}
\psi^{\#}(x|y)\chi^{\#}(y|x)\exp\{-F(y|y)/\hb\}\,dm(y).
$$
But the unity element of this multiplication is the section
$\cK$. To make this operation to be defined on functions, 
and to make the unity element to be the unity function~$1$, 
one 
must normalize this product in the following way:
$\psi*\chi\od(\cK\psi\times\cK\chi)/\cK$.
An explicit formula for the final product is \eqref{5}, 
that is, 
\begin{equation}
(\psi*\chi)(x)=\frac1{(2\pi\hb)^n}\int_{\cX}
p(x,y)\psi^{\#}(x|y)\chi^{\#}(y|x)\,dm(y).
\label{3.8}
\end{equation}

Another version of \eqref{3.8} is the formula for the left
multiplication 
\begin{equation}
\psi*=\frac1{\cK}\cdot\wh{\psi}\cdot\cK.
\label{3.9}
\end{equation}

The probability function is the ``eigenfunction'' of the
multiplication operator:
$$
\psi(x)*p(x,y)=p(x,y)*\psi(y)=\psi^{\#}(x|y)\,p(x,y),
$$
and, in particular, the following reproducing property holds:
$$
\psi(x)=\frac1{(2\pi\hb)^n}\int_{\cX}p(x,y)*\psi(y)\,dm(y),\qquad
\forall x\in\cX.
$$

Note that if $\psi\in L^1(\cX,dm)$ then the operator $\wh{\psi}$ is of
trace class: 
\begin{equation}
\tr\wh{\psi}=\frac1{(2\pi\hb)^n}\int_{\cX}\psi\,dm\,\od\,\tr\psi.
\label{3.10}
\end{equation}
{}For the product of two Wick operators we have
\begin{equation}
\tr(\wh{\psi}\,\wh{\chi}^*)=(\psi,\chi)\zs{W}
\label{3.11}
\end{equation}
and so the function space $W(\cX)$ is {\it isomorphic to the
space of Hilbert--Schmidt operators in\/} $\cL(\cX)$.
In particular, from \eqref{3.7} it follows
$$
|\psi(x)|\leq \|\psi\|\zs{W},\qquad \forall x\in \cX,
$$
and form \eqref{3.10}, \eqref{3.11} we have
\begin{equation}
\tr(\ochi*\psi)=\tr(\psi*\ochi)=(\psi,\chi)\zs{W},
\label{3.12}
\end{equation}
where the trace is defined  by the integral \eqref{3.10}. 

\begin{theorem}
The space $W(\cX)$ endowed with the product $*$ \eqref{3.8} is
an associative algebra with unity element~$1$ and with
involution $\psi\to\opsi$. 
Because of \eqref{3.12} $W(\cX)$ is a Frobenius algebra,
and moreover{\rm:}
$$
|(\psi*\chi)(x)|\leq \|\psi*\chi\|\zs{W}\leq 
\|\psi\|\zs{W}\cdot\|\chi\|\zs{W},\qquad \forall x\in\cX.
$$
The mapping $\psi\to\wh{\psi}$ given by \eqref{3.6} is an 
isomorphism of this algebra to the algebra of Hilbert--Schmidt
operators in $\cL(\cX)$, 
i.e., 
\begin{equation}
\wh{\psi}\cdot\wh{\chi}=\wh{\psi*\chi}.
\label{3.13}
\end{equation}
\end{theorem}

As usual, one can extend the algebra $W(\cX)$ in order to
include not only the Hilbert--Schmidt operators. 
Anyway, the product $*$ in any extended algebra $\cF_*(\cX)$ is 
given by the same formula \eqref{3.8} until the integral in
\eqref{3.8} makes sense (may be, as a distribution), 
the same is about the representation of this algebra by Wick
pseudodifferential operators given by~\eqref{3.6}.

Note that the Wick product \eqref{3.8} is represented by the
probability operator $\cP$ \eqref{3.1}:
\begin{equation}
(\psi*\chi)(x)=\cP_{y\to x}\big(\psi^{\#}(x|y)
\chi^{\#}(y|x)\big),
\label{3.14}
\end{equation}
where the subscript $y\to x$ indicates that the operator $\cP$
acts by the variable~$y$ and the result is a functions of the
variable~$x$.

If both the functions $\psi,\,\chi$ are smooth and do not depend on
the deformation parameter $\hb$ 
(or depend smoothly as $\hb\to0$)
and the reproducing measure also has the regular $\hb$-expansion
\eqref{2.18},  
then it is possible to derive from \eqref{3.14} the formal
$\hb$-power series expansion for the product $\psi*\chi$.
Indeed, 
when solving equation \eqref{2.8} for the reproducing
measure \eqref{2.11}, we have already obtained the expression:
\begin{equation}
\cP\sim I+\sum_{k\geq1}\hb^k\cP^{(k)},\qquad 
\cP^{(k)}=\sum^{k}_{s=0}f^{(k-s)}(x,\nabla)\circ \sigma_s,
\label{3.15}
\end{equation}
where the symbols $f^{(k)}$ and functions $\sigma_s$ are given
by \eqref{2.16}, \eqref{2.17} (we denote $f^{(0)}\equiv1$, 
$\sigma_0\equiv1$).
The coefficients $\cP^{(k)}$ of this expansion are differential
operators of order~$2k$ determined by the K\"ahler form $\omega$
only.  
{}For instance, 
\begin{equation}
\cP^{(1)}=\frac12\Delta\equiv \omega^{-1\mu\nu} \pa_\mu\opa_\nu,
\qquad
\cP^{(2)}=\frac18\Delta^2+\frac12\rho^{\mu\nu} \pa_\mu\opa_\nu,
\label{3.16}
\end{equation}
where $(\!(\rho^{\mu\nu})\!)$ is the Ricci tensor.
After the substitution of \eqref{3.15} to \eqref{3.14}, we obtain
the following formal $\hb$-power expansion for the Wick product:
\begin{equation}
(\psi*\chi)(x)\sim\psi(x)\chi(x)+\sum_{k\geq1}\hb^k\psi(x)
\underset{\leftrightarrow}{\cP}^{(k)}
\chi(x),\quad x\in\cX,
\label{3.17}
\end{equation}
where the differential operators 
$\cP^{(k)}=\underset{\leftrightarrow}{\cP}^{(k)}$ 
are defined by \eqref{3.15}; they act to the
left by holomorphic coordinates $z(x)$ and act to the right by 
antiholomorphic coordinates $\oz(x)$. 
A different way of calculation of the coefficients $\cP^{(k)}$
in \eqref{3.17}, based on the quantum tensor calculus, was 
suggested in [26]; other interesting approaches and a list of
references can be found in [10, 42].

{}Finally, we remark that, in general, the choice of the
reproducing measure~$dm$ is not unique. One can replace it by
the measure $(1+\psi_\hb)\,dm$, where $\psi_\hb$ is any function
from the null space $\Ker\cP$ such that $|\psi_\hb|<1$. 
This null space is certainly very big for the compact $\cX$,
but all such functions $\psi_\hb$ are highly oscillating
as $\hb\to0$ (see Example~4.1 below). 
So the assumption about 
the regular dependence of $dm$ as $\hb\to0$ is not redundant for
the expansion \eqref{3.17}.

\setcounter{equation}{0}

\section{Quantum restriction onto coadjoint orbits}

{}First of all, we consider the relationship between Wick
pseudodifferential operators and operators given by the
{\it geometric quantization theory}.

Let us introduce  the Poisson subalgebra $\cF^{(1)}(\cX)$
consisting of smooth functions over $\cX$ whose Hamiltonian flow
preserves the complex polarization. We denote by $\ad(\psi)$ the
Hamiltonian field corresponding to the function $\psi$ and split
this field in components along the complex polarization 
and along the conjugate one:
$$
\ad(\psi)=\ad_+(\psi)+\ad_-(\psi).
$$
If $\psi\in\cF^{(1)}(\cX)$ then the complex vector field
$\ad_-(\psi)$ transfers any antiholomorphic section to an
antiholomorphic section.

\begin{lemma}
If $\psi\in\cF^{(1)}(\cX)$ then the Wick pseudodifferential
operator $\wh{\psi}$ \eqref{3.6} is the following first order
differential operator{\rm:}  
\begin{equation}
\wh{\psi}=\psi+i\ad_-(\psi)(F)-i\hb\ad_-(\psi).
\label{4.1}
\end{equation}
Moreover, the Dirac axiom holds{\rm:}
\begin{equation}
\frac{i}{\hb}[\wh{\psi},\wh{\chi}]=\wh{\vbox to 12pt{}\{\psi,\chi\}},
\qquad \forall \psi,\chi\in\cF^{(1)}(\cX).
\label{4.2}
\end{equation}
\end{lemma}

Note that \eqref{4.1} is exactly the construction of the
geometric quantization over K\"ahlerian manifolds (see [34]).
So, the Wick pseudodifferential calculus extends the geometric
quantization scheme from the Lie algebra level \eqref{4.2} to
the associative algebra level \eqref{3.13}. 
In particular, for generic $\psi\not\in\cF^{(1)}(\cX)$, 
the operator $\wh{\psi}$ ceases to be the first order
differential operator. 

\begin{example}
\rm
Let $\cX$ be the unit sphere 
$$
\bS^2=\{\xi\cdot \xi=1\}\subset\bR^3,\quad 
\xi=(\xi^1,\xi^2,\xi^3), \,\,\text{$\xi^j$ are Cartesian coordinates}.
$$
The complex structure is standard:
$z=\frac{\xi^1+i\xi^2}{1-\xi^3}$,
$\omega=\omega_0=\frac{2id\oz\wedge dz}{(1+|z|^2)^2}$.
The quantization condition \eqref{2.2} implies 
$\hb=2/N$, where $N\in\bZ_+$. 
The reproducing measure is $dm=(1+\hb/2)\,dm^\omega$, and the
probability function is  
$p=\big(\frac{1+\xi\cdot \xi'}{2}\big)^N$,
where $\xi,\xi'\in\bS^2$.
The probability operator is given by the formula 
$$
\cP=\sum^{N}_{k=0}\frac{(N+1)!N!}{(N+1+k)!(N-k)!}P^{(k)},
$$
where $P^{(k)}$ is the projection to the space $\cF^{(k)}$ of
$k$th spherical harmonics, that is to 
the eigenspace of the Laplacian $\Delta_{\bS^2}$ corresponding
to the eigenvalue $k(k+\nobreak1)$.
In this case 
$$
\Ker\cP=\bigoplus_{k\geq N+1}\cF^{(k)},\qquad 
\cF_*(\bS^2)=\bigoplus_{0\leq k\leq N}\cF^{(k)},\qquad 
\dim \cF_*(\bS^2)=(N+1)^2.
$$
{}For instance, if $N=1$ (that is, $\hb=2$), 
then $\dim\cF_*(\bS^2)=4$.
In this case let us take the following basis in $\cF_*(\bS^2)$:
$1$, $x^j\od \xi^j\big|_{\bS^2}$ \ 
$(j=1,2,3)$.
Then the Wick algebra structure \eqref{3.8} on $\cF_*(\bS^2)$ is
generated by relations between basis functions:
\begin{align*}
x^j*x^j&=1 \qquad (j=1,2,3),\\
x^1*x^2&=-x^2*x^1=-i x^3,\\
x^2*x^3&=-x^3*x^2=-i x^1,\\
x^3*x^1&=-x^1*x^3=-i x^2.
\end{align*}
So, in this case $\cF_*(\bS^2)$ is {\it isomorphic to the quaternion
algebra}.  

Note that for arbitrary $N$ all functions $x^j$ belong to the
Poisson subalgebra $\cF^{(1)}(\bS^2)\approx \su(2)$, and the
first order differential operators $\wh{x}^j$ \eqref{4.1} 
represent this Lie algebra in the $(N+1)$-dimensional Hilbert
space $\cL(\bS^2)$: 
$$
\frac{i}{\hb}[\wh{x}^1,\wh{x}^2]=\wh{x}^3,\quad
\frac{i}{\hb}[\wh{x}^2,\wh{x}^3]=\wh{x}^1,\quad
\frac{i}{\hb}[\wh{x}^3,\wh{x}^1]=\wh{x}^2,\quad
\wh{x}^{j\,*}=\wh{x}^j\quad (j=1,2,3).
$$
The explicit formulas for $\wh{x}^j$ are the following:
$$
\wh{x}^1-i\wh{x}^2=\hb\oz(\oz\opa-N),\qquad 
\wh{x}^1+i\wh{x}^2=-\hb\opa,\qquad 
\wh{x}^3=\hb(\oz\opa-N/2).
$$
The direct derivation of the Casimir operator gives
\begin{equation}
(\wh{x}^1)^2+(\wh{x}^2)^2+(\wh{x}^3)^2
=(1+\hb)\cdot I\ne I.
\label{4.3}
\end{equation}
But the classical Casimir is 
$(x^1)^2+(x^2)^2+(x^3)^2=\xi\cdot \xi\big|_{\bS^2}=1$. The
difference between classical and quantum values of the Casimir
elements is one of the reasons which gives rise to the question
about the quantum restriction operation~\eqref{3}.
\end{example}

Now let us consider a general Lie
algebra $\fg$ and its coadjoint orbits [31] $\cX\subset \fg^*$
with standard symplectic form $\omega=\omega_0$ generated on
$\cX$ by the linear Lie--Poisson brackets from $\fg^*$. 
We assume that $\dim\cX$ is maximal and  there is a
$\fg$-invariant complex structure on $\cX$ with respect to which the
form $\omega_0$ is K\"ahlerian
(for instance, it is true if $\fg$ is compact). 
Assume that the quantization
condition \eqref{2.2} holds 
and construct the space $\cL(\cX)$ using K\"ahler potentials 
$F=F_0$ of the form $\omega_0$ and the measure 
$dm=\const\cdot dm^{\omega_0}$.  

Any vector $\eta\in\fg$ is identified with the linear function over
$\fg^*$ by the formula 
$\eta(\xi)\od\langle \eta,\xi\rangle$ \ $\forall\xi\in\fg^*$. 
Consider the classical restriction $\eta\big|\zz{\cX}$ of the
function $\eta$ to the orbit $\cX$.  
Then $\eta\big|\zz{\cX}\in\cF^{(1)}(\cX)$, and one can apply the
construction of the geometric quantization from Lemma~4.1. So,
the first order differential operator 
$\wh{\eta\big|\zz{\cX}}$ appears in the Hilbert space $\cL(\cX)$ of
antiholomorphic sections over~$\cX$.

Let us fix any basis $\eta^1,\dots,\eta^m$ in $\fg$ and denote 
$x^j=\eta^j\big|\zz{\cX}$. Then we have 
the set of operators $\wh{x}^j$ in $\cL(\cX)$.

Each polynomial $f$ on $\fg^*\approx\bR^m$ is given by a
polynomial function of the basis elements 
$f=f(\eta^1,\dots,\eta^m)$. 
Its classical restriction onto the coadjoint orbit is given by 
$f\big|\zz{\cX}=f(x^1,\dots,x^m)$. 

The {\it quantum restriction\/} $f\big|\zz{\wh{\cX}}$ we define 
in such a way that  
\begin{equation}
\wh{f\big|\zs{\wh{\cX}}}=f(\wh{x}^1,\dots,\wh{x}^{\,m}),
\label{4.4}
\end{equation}
where on the right-hand side the operators $\wh{x}^j$ are
Weyl-symmetrized. 

The quantum product of two Weyl-symmetrized polynomials on
$\fg^*$ we denote by $\star$. So, in view of \eqref{4.4}, we have
$$
\wh{f\big|\zs{\wh{\cX}}}\cdot\wh{g\big|\zs{\wh{\cX}}}
=\wh{(f\star g)\big|\zs{\wh{\cX}}}.
$$
Using the Wick product $*$ of symbols over the K\"ahlerian
manifolds $\cX$, we obtain the desirable formula \eqref{4},
i.e., 
the homomorphism $f\to f\big|\zz{\wh{\cX}} $ from the algebra of
polynomials $\cF_\star(\fg^*)$ to the Wick algebra 
$\cF_*(\cX)$.

To calculate the quantum restriction $f\big|\zz{\wh{\cX}}$
explicitly,  
let us apply formula \eqref{3.7} to the relation \eqref{4.4}: 
$$
f\big|\zs{\wh{\cX}}=\tr\big(f(\wh{x}^1,\dots,\wh{x}^{\,m})\Pi\,\big)
=\frac1{\cK}f(\wh{x}^1,\dots,\wh{x}^{\,m})\cK.
$$
Since $\cK=\exp(F_0/\hb)$ and $\wh{x}^j$ are defined by
\eqref{4.1} with $F=F_0$, we obtain the 
following statement.

\begin{theorem}
The quantum restriction of the function $f\in\cF(\fg^*)$ onto
the coadjoint orbit $\cX$ with invariant complex structure
is given by the formula
\begin{equation}
f\big|\zs{\wh{\cX}}=f\big(x-i\hb\ad_-(x)\big)\,1.
\label{4.5}
\end{equation}
On the right-hand side of \eqref{4.5} the first order
differential operators, which are arguments of the function~$f$,
are Weyl symmetrized, and the resulting operator is applied to
the unity function~$1$ over~$\cX$.
\end{theorem}

Let us calculate the leading quantum correction of order $\hb$ 
in the formal $\hb$-expansion \eqref{3} for the quantum
restriction. 
{}From \eqref{4.5} it follows that
\begin{equation}
e_1=\frac{1}{2} \bg^{j\ell}_-D_jD_\ell
=\frac14(\Delta-\Delta(x^j)D_j).
\label{4.6}
\end{equation}
Here $D_j$ are partial derivatives by the coordinates 
$\xi^j=\langle\eta^j,\xi\rangle$ 
in $\fg^*$; by $\Delta$ we denote the Laplace--Beltrami operator on the
K\"ahler orbit $\cX\subset \fg^*$; and
$\bg^{j\ell}_-$ is the ``polarized part'' of the
Lie--Poisson tensor, i.e.,  
\begin{equation}
\bg^{j\ell}_-\,\od\,\omega^{-1\mu\nu}_0\pa_\mu x^j \cdot \opa_\nu x^\ell
=-i\ad_-(x^j)(x^\ell).
\label{4.7}
\end{equation}
The tensor \eqref{4.7} is Hermitian and degenerate,
$\rank\bg_-=\dim\cX$. 

{}For instance, if $f$ is a {\it quadratic function\/} 
on~$\fg^*$, then its quantum restriction to the orbit~$\cX$
is precisely given by 
$$
f\big|\zs{\wh{\cX}}=f\big|\zs{\cX}+\frac{\hb}{2}\tr(\bg_-D^2f).
$$

\setcounter{example}{0}

\begin{example} \rm (continuation).
In the case $\fg=\su(2)$, $\cX=\bS^2$ we have
\begin{equation}
e_1=\frac14\Delta_{\bS^2}+\frac{\pa}{\pa K},\qquad 
K\od (\xi^1)^2+(\xi^2)^2+(\xi^3)^2.
\label{4.8}
\end{equation}
Since $K$ itself is a quadratic function, we obtain the precise
formula:  
$$
K\big|\zs{\wh{\cX}}=K\big|\zs{\cX}+\hb e_1(K)=1+\hb.
$$
This result exactly correlates with the derivation
\eqref{4.3}.
\end{example}

We denote by $\{K_s\}$ a basis of Casimir functions of Lie--Poisson
algebra $\cF(\fg^*)$; 
then coadjoint orbits are determined by equations 
$\{K_s=\const\}$, and $\{\pa/\pa K_s\}$ is a basis of normal
vector fields over each orbit of maximal dimension. 
Then we can 
split the operator $e_1$ \eqref{4.6} in components tangent and 
transversal to~$\cX$:
$e_1=e^{\|}_1+e^{\perp}_1$,
where the transversal component is the vector field given by  
\begin{equation}
e^{\perp}_1=-\frac14\Delta(x^j)\cdot D_j K_s\frac{\pa}{\pa K_s}
=\frac12  \bg^{j\ell}_-
D^2_{j\ell} K_s \frac{\pa}{\pa K_s}.
\label{4.9}
\end{equation}

Note that the expression on the right-hand side of \eqref{4.9}
does not depend on the choice of Casimir functions~$K_s$.

Thus for compact Lie algebras and their coadjoint orbits of
maximal dimension we obtain the following statement.

\begin{corollary}
The leading quantum correction $e_1$ in the 
$\hb$-expansion \eqref{3} of the quantum restriction onto the
coadjoint orbit $\cX\subset \fg^*$ is the 
second order differential operator given by \eqref{4.6}.
The component normal to the orbit is 
the first order operator given by \eqref{4.8}. 
This normal vector field $e^\perp_1$ is correctly defined 
on the domain $\cN$ of regular points in $\fg^*$, 
where the rank of the Lie-Poisson tensor is maximal.

The field $e^\perp_1$ differs only  by a constant multiplier
from  the mean curvature vector field determined by the 
Levi--Civita connection on K\"ahler leaves $\cX\subset \fg^*$
and by the affine connection on $\fg^*$. 
\end{corollary}

The last statement of this Corollary was conjectured 
by A.~Weinstein in the discussion about formulas \eqref{4.6},
\eqref{4.8} presented in the author's lecture at the conference 
``Poisson--2000'' (Luminy, June 2000). 

{}From the viewpoint of the Hochschild complex [2, 17],
the quantum corrections $e_k$ in~(1.3) are $1$-cochains 
obeying the series of equations 
$\delta e_k=c^*_k-c^\star_k-\mu_k$. 
Here $\delta$ denotes the Hochschild differential, 
and $c_k$ are $2$-cochains from expansion~(1.1);
upper signs $*$ and $\star$ labling the algebra: 
$\cF_*(\cX)$ and $\cF_\star(\cU)$.
The $2$-cochains $\mu_k$ are determined by previous 
$c_1,\dots,c_{k-1}$ and $e_1,\dots,e_{k-1}$;
for instance, $\mu_1=0$, 
$\mu_2(f,g)=e_1(c^\star_1(f,g))
-c^*_1(f,e_1(g))-c^*_1(e_1(f),g)$.
They obey the equations 
$\delta\mu_k=\nu^*_k-\nu^\star_k$, where $\nu_k=\delta c_k$. 
The $3$-cocycles $\nu_k$ are given by the previous 
$c_1,\dots,c_{k-1}$; for instance, $\nu_1=0$, 
$\nu_2(f,g,k)=c_1(f,c_1(g,k))-c_1(e_1(f,g),k)$. 

Now we calculate the quantum restriction \eqref{4.5} 
without any $\hb$-expansions by using expression (1.3a).
The symbol $E\zs{\cX}$ of the operator of quantum restriction is
given by
\begin{equation}
E\zs{\cX}(\,\cdot,\eta)
=e^{-i\eta/\hb}\big|\zs{\cX}\cdot e^{i\eta/\hb}\big|\zs{\wh{\cX}}.
\label{4.10}
\end{equation}
On the right-hand side we consider the vector $\eta\in\fg$ as a
linear function on $\fg^*$ and apply both types of restriction
onto the orbit~$\cX$: classical and quantum.
{}For the quantum one let us apply Theorem~4.1, that is, take 
$f=\exp(i\eta/\hb)$ in \eqref{4.5}:
$$
e^{i\eta/\hb}\Big|\zs{\wh{\cX}}=\exp\Big\{\frac i\hb
\langle\eta,\,x-i\hb\ad_-(x)\rangle\Big\}1.
$$
The operators $\ad_-(x)$ are of the first order, and so the
action of the last exponent can be evaluated by the method of
characteristics.  
Thus we obtain from\eqref{4.10}:
\begin{equation}
E\zs{\cX}(\xi,\eta)=\exp\bigg\{\frac i\hb
\bigg(\int^1_0\langle\eta,\Xi_-(t,\xi,\eta)\rangle\,dt 
- \langle\eta,\xi\rangle\bigg)\bigg\}
=\exp\bigg\{-\frac{\eta\wt{\bg}\eta}\hb\bigg\}.
\label{4.11}
\end{equation}
Here $\Xi_-$ is the characteristic:
\begin{equation}
\frac{d}{dt}\Xi_-=i\eta \bg_-(\Xi_-),\qquad 
\Xi_-\Big|_{t=0}=\xi\in\cX,
\label{4.12}
\end{equation}
the tensor $\bg_-$ is defined by \eqref{4.7},
and 
$$
\wt{\bg}(\xi,\eta)\,\od\,\int^1_0(1-t)\bg_-\big(\Xi_-(t,\xi,\eta)\big)\,dt.
$$

\begin{corollary}
{\rm(a)} The quantum restriction onto the coadjoint orbit 
$\cX$ is given by formula {\rm(1.3a)}, 
where the symbol $E\zz{\cX}$ is determined by \eqref{4.11}
via the trajectories $\Xi_-$ of the ``polarized''
Lie--Hamilton system \eqref{4.12}. 

{\rm(b)} The group elements \eqref{12} are given by the formula
\begin{equation}
e^{i\eta/\hb}_*(\xi)=
\exp\bigg\{-\frac 1\hb\eta\wt{\bg}(\xi,\eta)\eta
+\frac i\hb\langle\eta,\xi\rangle\bigg\},\qquad 
\xi\in\cX.
\label{4.13}
\end{equation}
\end{corollary}

Some other way to calculate the exponential function from
Corollary 4.2\,(b) was presented in [24, 25] using the areas of
dynamical membranes in~$\cX^{\#}$.

\hbox to 1cm{}

Note that the trajectories $\Xi_-$ in \eqref{4.12} leave the
real coadjoint orbit $\cX$ and belong to its complexification. 
Nevertheless, the action functions in the exponents
\eqref{4.11}, \eqref{4.13}
have a nonnegative imaginary part. The group elements 
$e^{i\eta/\hb}_*$ are asymptotically (as $\hb\to0$) concentrated
at the points where this imaginary part vanishes, namely, 
at fixed points of the coadjoint action of $\exp(\eta)$
on~$\cX$. This {\it effect of concentration\/} is purely quantum
one.  
The quantum group element $e^{i\eta/\hb}_*$ is exponentially small as
$\hb\to0$ outside the set of fixed points. 
On the other hand, 
for fixed $\xi\in\cX$, the element (4.13) is 
{\it an exact eigenfunction of the Laplace operators\/} 
on the Lie group $G\sim \exp(\fg)$. The oscillation front
of this function as $\hbar\to0$ is an isotopic 
submanifold in $T^*G$ generated by the stabilizer $G_\xi\subset G$.


\setcounter{example}{0}

\begin{example} \rm (continuation).
In the case $\fg=\su(2)$, $\cX=\bS^2=\{\xi\cdot\xi=1\}$ we have 
$$
\bg^{j\ell}_-(\xi)=\frac12(\delta^{j\ell}-\xi^j\xi^\ell).
$$
Instead of the vector characteristic $\Xi_-$ we consider the
scalar $X=\langle\eta,\Xi_-\rangle$. Then equations \eqref{4.12} are
reduced to  
$$
\frac{d}{dt}X=\frac{i}{2}\big(|\eta|^2-X^2\big),\qquad
X\Big|_{t=0}=\langle\eta,\xi\rangle.
$$
The solution is
$$
X=|\eta|\frac{(\langle\eta,\xi\rangle+|\eta|)\exp\{i|\eta|t\}
+\langle\eta,\xi\rangle-|\eta|}
{(\langle\eta,\xi\rangle+|\eta|)\exp\{i|\eta|t\}
-\langle\eta,\xi\rangle+|\eta|}.
$$
Thus from formula \eqref{4.11} we obtain the symbol
$E\zz{\bS^2}$ of the quantum restriction 
operation (1.3a):
\begin{align}
E\zz{\bS^2}(\xi,\eta)
&=
\exp\bigg\{\frac i\hb\bigg(\int^1_0X(t,\xi,\eta)\,dt
-\langle\eta,\xi\rangle\bigg)\bigg\}\nonumber\\
&=\bigg(\cos\frac{|\eta|}{2}+\frac{i}{|\eta|}\sin\frac{|\eta|}{2}
\cdot\langle\eta,\xi\rangle\bigg)^{2/\hb}
e^{-i\langle\eta,\xi\rangle/\hb}
\label{4.14}
\end{align}
(we recall that $2/\hb=N$ is an integer number, and
$\xi\in\bS^2$).  
The approximation of $E\zz{\bS^2}(\xi,\eta)$ near $\eta=0$ is
the following: 
\begin{equation}
E\zz{\bS^2}(\xi,\eta)=\exp\bigg\{
-\frac1{4\hb}\big(|\eta|^2-\langle\eta,\xi\rangle^2\big)
+O(\eta^3/\hb)\bigg\}.
\label{4.15}
\end{equation}
So, the quantum restriction operator 
$f\big|\zz{\wh{\bS^2}}=\Big(E\zz{\bS^2}(\stackrel{2}{\xi},\,
-i\hb\stackrel{1}{\pa/\pa\xi})f(\xi)\Big)\big|\zz{\bS^2}$ is given
precisely by formula \eqref{4.14}, and its first 
approximation  (on $\hb$-in\-de\-pen\-dent functions)
is evaluated by \eqref{4.15}:
$$
E\zz{\bS^2}(\stackrel{2}{\xi},\,-i\hb\stackrel{1}{\pa/\pa\xi})
=I+\frac{\hb}{4}\sum^{3}_{s,\ell=1}(\delta^{s\ell}-\xi^s\xi^\ell)
\frac{\pa^2}{\pa\xi^s\pa\xi^\ell}+O(\hb^2).
$$
{}From this formula we obtain, of course, the same expression
\eqref{4.9} for the first quantum correction~$e_1$.
All other corrections $e_j$ ($j\geq2$) in expansion \eqref{3}
are easily extracted from the precise formula \eqref{4.14}.

The group elements in this example are derived from
\eqref{4.14}: 
\begin{equation}
e^{i\eta/\hb}_*=E\zs{\bS^2}e^{i\eta/\hb}
=\bigg(\frac{\cos(|\eta|/2)}{\cos(S_\eta/2)}\bigg)^{2/\hb}
e^{iS_\eta/\hb},
\label{4.16}
\end{equation}
where $S_\eta(\xi)=2\arctan(\langle\frac{\eta}{|\eta|},\xi\rangle
|\tan \frac{|\eta|}{2}|)$ and $2/\hb=N$. 
This formula presents a realization of the Lie group $G=\SU(2)$ 
in the function algebra $\cF_*(\bS^2)$ 
(of spherical harmonics of order $\leq N$).  
The corresponding Wick differential operators make up 
the Hermitian irreducible representation $\pi\zz{\bS^2}$ of the
group $\SU(2)$ 
in the space $\cL(\bS^2)$ (of all polynomials of degree $\leq N$). 

Let the vector $\eta\ne0$ belong to the domain 
$\{|\eta|<2\pi\}$ where the exponential mapping 
$\exp:\,\su(2)\to \SU(2)$ is one-to-one. 
Then the function \eqref{4.16} is concentrated at two points 
$\xi_{\pm}=\pm\eta/|\eta|$,
and it is exponentially small (as $\hb\to0$) out of this set.
The points $\xi_{\pm}$ are fixed points of the coadjoint action
of the element $\exp(\eta)\in\SU(2)$ on the orbit $\cX=\bS^2$.
Of course, the exceptional value $\eta=0$ 
corresponds to the unity element of the group:
in this case the function $e^{i\eta/\hb}_*$ is equal to~$1$
identically. 
On the boundary circle $|\eta|=2\pi$ 
the function $e^{i\eta/\hb}_*$ is equal to $(-1)^N$, and so, 
if~$N$ is even (what corresponds to 
representations of $\SO(3)$) then $e^{i\eta/\hb}_*$ is the same
unity function~$1$. But if~$N$ is odd then one has to go to the
second sheet of the universal covering, that is, to the domain 
$2\pi<|\eta|<4\pi$, in order to obtain the realization of the
whole group $\SU(2)$ in the function space $\cF_*(\bS^2)$.  
\end{example}

\setcounter{equation}{0}

\section{Irreducible representations and special\\ functions
corresponding to complex\\ polarizations}

As we saw in Section~4, the construction of quantum restriction
to coadjoint orbits is actually equivalent to the construction
of irreducible representations corresponding to these orbits 
(formulas \eqref{8}, \eqref{4.4} demonstrate this relationship).
Until we are in the case of Lie algebras, we can use the
geometric quantization theory to produce these irreducible
representations. 
But if we deal with nonlinear Poisson brackets and with algebras
$\cF_\star(\cM)$ generated by non-Lie permutation relations, 
then the problem of quantum restriction onto symplectic leaves
$\cX\subset\cM$ becomes rather difficult, since in this case
there is no general construction of irreducible representations. 

In the given section, following [26], we describe one 
possible construction based on the notion of complex
polarization of quantum algebra. 
In what follows, we modify the approach [26] 
in order to avoid the use of the complexified phase space 
(complexified symplectic groupoid) over $\cM$.

Let $\cF=\cF(\cM)$ be a quantized algebra of functions over
$\cM$ with an associative multiplication $\star$ satisfying
\eqref{1}, with the unity element, and the involution $f\to\of$.
{}For any two subspaces $\cR$, $\cT\subset\cF$ 
we denote by $\cR\ostar\cT$ the subspace
otained by the composition of tensor product and multiplication
~$\star$. 

Let $\cF^+$ be a complex polarization in $\cF$, i.e., 
a subalgebra with the following properties:
\begin{description}
\item[(i)] the subspace $\overline{\cF^+}\ostar\cF^+$ 
coincides with~$\cF$, 
\item[(ii)] the subalgebra
$\cF^0\,\od\,\cF^+\cap\overline{\cF^+}$ is commutative.
\end{description}
A point $a\in\cM$ is called a {\it vacuum point\/} with respect
to the polarization if  
\begin{description}
\item[(iii)] 
$(f\star g)(a)=f(a)g(a)$ \ 
$\forall\,f\in\cF$, \ $\forall\,g\in\cF^+$.
\end{description}
Under these conditions we call $\star$ the {\it normal product}. 

\hbox to 1cm{}

Let us choose a commutative functional basis $Z^1,\dots,Z^n$ 
in $\cF^+\setminus\cF^0$
and denote by $P_1,\dots, P_n$ the Darboux dual subset 
of functions in~$\cF$:
\begin{equation}
[Z^s,P_j]\zz\star=\hb\delta^s_j,\qquad 
[Z^s,Z^j]\zz\star=[P_s,P_j]\zz\star=0\quad (s,j=1,\dots,n).
\label{5.1}
\end{equation}
To simplify the construction, we assume that $Z^s$ and $P_j$
satisfy the boundary conditions 
\begin{equation}
Z^s(a)=0,\qquad P_j(a)=0\quad (s,j=1,\dots,n).
\label{5.2}
\end{equation}

Also assume that the following $\star$-exponential exists in the
algebra $\cF$: 
\begin{equation}
U_z\,\od\,\exp_\star(zP/\hb),\qquad z=(z^1,\dots z^n)\in\cD.
\label{5.3}
\end{equation}
Here $\cD\subset\bC^n$ is a domain containing the point $z=0$. 
The function $U_z$ is the solution of the
Cauchy problem  
\begin{equation}
\hb\pa U_z/\pa z=P\star U_z,\qquad U_0=1
\label{5.4}
\end{equation}
(the general $\star$-exponential was considered in [3]; see the
discussion in [19]).
We denote
\begin{equation}
\cK(\oz|z)\,\od\,(\oU_z\star U_z)(a)
=\sum_{|\alpha|,|\beta|\geq0}
\frac{\oz^\alpha z^\beta}{\hb^{|\alpha|+|\beta|}}
k_{\alpha\beta},
\label{5.5}
\end{equation}
where
$$
k_{\alpha\beta}=\frac1{\alpha!\beta!}
(\oP^{\star\alpha}\star P^{\star\beta})(a).
$$
The sign $\star$ in notation of powers means that these powers
are taken in the algebra $\cF$ (i.e., in the sense of the normal
product $\star$).
In addition to the above properties (i)--(iii),
we introduce the quantum K\"ahlerian condition:
\begin{description}
\item[(iv)] the matrix $(\!(k_{\alpha\beta})\!)$ is positive
definite (maybe, not strictly). 
\end{description}

By $\cL_a$ we denote the space of antiholomorphic distributions  
$u(\oz)=\sum{\oz^\alpha}u_\alpha$
generated by vectors $(\!(u_\alpha)\!)$ orthogonal to the null
kernel of the matrix $(\!(k_{\alpha\beta})\!)$, and on $\cL_a$
introduce the Hilbert norm 
\begin{equation}
\|u\|\,\od\,\Big(\sum\hb^{|\alpha|+|\beta|}k^{-1\alpha\beta}
u_\alpha\ou_\beta\Big)^{1/2}.
\label{5.6}
\end{equation}

Obviously, the function $\cK$ \eqref{5.5} is just the
reproducing kernel of the space $\cL_a$, that is, the integral
kernel of the unitary operator in $\cL_a$.
One can ask about the reproducing measure $dm$ for
the space $\cL_a$:
\begin{equation}
\|u\|=\bigg(\frac1{(2\pi\hb)^n}\int\frac{|u|^2}{\cK}\,dm\bigg)^{1/2},
\label{5.7}
\end{equation}
where the integral is taken over the domain~$\cD$. 
Then, in the space $\cL_a$ we can define the Wick
pseudodifferential operators using the kernel $\cK$ and the
Hilbert structure \eqref{5.6} or \eqref{5.7}. 
We denote these operators again by the hat sign: $\wh{\psi}$, 
where $\psi$ are symbols on $\cD$ extended to symbols 
$\psi^{\#}(\oz|z)$ on $\cD\times\cD$ holomorphic in~$z$ and
in~$\oz$. 

Let us take $f\in\cF(\cM)$ and determine the following Wick
symbol   
\begin{equation}
f_a(\oz|z)=\frac{(\oU_z\star f\star U_z)(a)}{(\oU_z\star U_z)(a)}.
\label{5.8}
\end{equation}

\begin{theorem}
The correspondence 
\begin{equation}
f\to\wh{f}_a
\label{5.9}
\end{equation}
is an irreducible Hermitian representation of the algebra
$\cF(\cM)$ {\rm(}with the normal product $\star${\rm)}
in the Hilbert space~$\cL_a$.
\end{theorem}

Explicit formulas for operators $\wh{f}_a$ were given in [26]
via symbols of left and right quantum reduction mappings [22,
27] in the 
complexified symplectic groupoid over $\cM$.

Consider what happens with coordinate functions $Z^s\in\cF^+$
under the correspondence \eqref{5.9}.
{}From \eqref{5.8} and \eqref{5.1} it follows that 
\begin{equation}
(Z^s)_a=z^s,\qquad s=1,\dots,n.
\label{5.10}
\end{equation}
Thus the representation \eqref{5.9} transforms each function 
$\oZ^s$ to the operator of multiplication by $\oz^s$, and each
function $Z^s$ to the adjoint operator $\wh{z}^{\,s}=(\oz^s)^*$ in
the Hilbert space $\cL_a$.

Now let us take a basis of functions $A^1,\dots,A^k$ in the
subalgebra $\cF^0$ (from condition (ii))
and pass to the classical limit:
$$
\cA^j(\oz|z)\,\od\,\lim_{\hb\to0}(A^j)_a(\oz|z),\qquad j=1,\dots,k.
$$
Explicit formulas for functions $\cA^j$ are given in [26].

Note that $Z^1,\dots,Z^n$ are complex coordinates, and
$A^1,\dots,A^k$ are real coordinates on $\cM$, 
and $2n+k=\dim\cM$, $2n=\dim\cX_a$.

\begin{corollary}
In the Poisson manifold $\cM$ the symplectic leaf $\cX_a$
{\rm(}containing the vacuum point~$a${\rm)} 
is given by the following parametrization of coordinates{\rm:}
\begin{equation}
Z=z,\qquad A=\cA(\oz|z).
\label{5.11}
\end{equation}
Here $z$ is running over a domain $\cD\subset\bC^n$.
{}Formula~\eqref{5.11} introduces the complex
structure to the leaf $\cX_a$ from the local chart~$\cD$, 
and transfers each symbol \eqref{5.8} onto $\cX_a$. 
The symbols $f_a$ are quantum restrictions of functions~$f$ to
the leaf{\rm:} 
$$
f\big|\zs{\wh{\cX}_a}=f_a,\qquad 
(f\star g)\big|\zs{\wh{\cX}_a}=
f\big|\zs{\wh{\cX}_a}*g\big|\zs{\wh{\cX}_a},
$$
where $*$ is the Wick product over $\cX_a$ generated by
the K\"ahlerian form 
\begin{equation}
\omega=i\hb\opa \pa \ln\cK.
\label{5.12}
\end{equation}
If the reproducing measure $dm$ \eqref{5.7} exists, then the
space $\cL_a$ is identified with the space $\cL(\cX_a)$ of
antiholomorphic sections.
The quantum K\"ahlerian objects on $\cX_a$ 
in the classical limit $\hb\to0 $ admit the asymptotics 
\begin{equation}
\omega=\omega_0+O(\hb),\qquad dm=dm^{\omega_0}+O(\hb),
\label{5.13}
\end{equation}
where the symplectic form $\omega_0$ on $\cX_a$ is generated
from $\cM$ by the Poisson structure.
\end{corollary}

The reproducing kernel $\cK$ \eqref{5.5} is the key object of
the construction just described.
This kernel determines the quantum K\"ahler structure
\eqref{5.12} on the symplectic leaf $\cX_a$ satisfying the 
quantization condition $\frac1{2\pi\hb}[\omega]\in 
H^2(\cX_a,\bZ)$.
Note that $\cK$, as a holomorphic section of the complex line 
bundle over $\cX^{\#}_a$, does not depend 
on the choice of the bases $\{Z^j\}$ and $\{P_j\}$
(modulo changes of variables and transfers to equivalent
bundles).  
That is why we call the kernel \eqref{5.5} a {\it special function\/}
corresponding to the quantum complex polarization over the
symplectic leaf $\cX_a$ with the vacuum point~$a$.

In the next section we demonstrate some examples where $\cK$
turns out to be a hypergeometric function.

\setcounter{equation}{0}

\section{Quantum surfaces of revolution and\\ hypergeometric functions}

As an example we consider two-dimensional leaves (surfaces) of
simplest  topology but with arbitrary complex structure.
Namely, let us consider the algebraic surface 
\begin{equation}
\rho(t)-(S^2_1+S^2_2)=K=\const,
\label{6.1}
\end{equation}
where $S_1,\,S_2,\,t$ are coordinates in $\bR^3$ and
the function $\rho$ is a polynomial.
{}For each $t$ Eq.~\eqref{6.1} describes a circle, so \eqref{6.1}
is a surface of 
revolution with the axis~$t$. 
The topology depends on values of the constant~$K$.
If $\max\rho>K>\min\rho$, then the surface is homeomorphic to the plane
or to the sphere, 
if $K<\min\rho$, then the surface is homeomorphic to the
cylinder or the torus.

It is more convenient to deal with surfaces embedded in $\bR^{k+2}$ 
(with $k\geq1$) by equations
\begin{align}
&\rho(A)-(S^2_1+S^2_2)=K=\const,\label{6.2}\\
&\varkappa_j(A)=\const_j\qquad (j=1,\dots,k-1),
\nonumber
\end{align}
where $A=(A_1,\dots,A_k), A_\mu$ are coordinates in $\bR^k$. 
Additional independent polynomial functions
$\varkappa_1,\dots,\varkappa_{k-1}$ in \eqref{6.2} 
determine an algebraic curve in $\bR^k$ playing the role of the
axis of revolution. Note that somewhere we need to consider not 
the whole $\bR^k$
but only a suitable domain in $\bR^k$; we shall do this without
making additional notations.

Our goal is to introduce a quantum version of the surfaces
\eqref{6.2}. {}First of all, one has to quantize the space 
$\bR^{k+2}$ to where these surfaces are embedded. 
To do this, let us introduce  a Poisson structure
to $\bR^{k+2}$ such that surfaces \eqref{6.2} are symplectic
leaves (i.e., $K$ and $\varkappa_j$ are Casimir functions).

Let $v$ be a vector field on $\bR^k$ annulling all $\varkappa_j$;
that is, $\varkappa_j$ are integrals of motion for the
dynamical flow 
\begin{equation}
\gamma^t:\,\,\bR^k\to\bR^k,\qquad
\frac{d}{dt}\gamma^t=v(\gamma^t),\qquad
\gamma^0=\id.
\label{6.3}
\end{equation}
If we like to live inside the algebraic case, we have to assume
that this flow preserves polynomials. It is certainly true if
all the components $v_\mu(A)$ of the field  
$v=v_\mu(A)\pa/\pa A_\mu$ are polynomial, and 
also the Jacobi matrix $(\!(\pa v_\mu/\pa A_\nu\,)\!)$ 
is nilpotent for all $A\in\bR^k$;
then the trajectories 
$\gamma^t(A)$ are polynomial in~$A$ and in~$t$.

Denote $\lambda=v(\rho)$, and introduce the following Poisson
brackets on $\bR^{k+2}$:
\begin{gather}
\{S_2,S_1\}=\frac12\lambda(A),\quad
\{A_\mu,S_1\}=v_\mu(A)S_2,\quad
\{A_\mu,S_2\}=-v_\mu(A)S_2,
\label{6.4}\\
\{A_\mu,A_\nu\}=0,\quad \mu,\nu=1,\dots,k.
\nonumber
\end{gather}

\begin{lemma}
Symplectic leaves of the Poisson structure \eqref{6.4} are given
by Eqs.~\eqref{6.2}.
\end{lemma}

If we denote $C=\oB=S_1+iS_2$, then the nontrivial relations
\eqref{6.4} take the form
$$
\{C,B\}=i\lambda(A),\qquad \{C,A\}=iv(A)C.
\eqno{(6.4\text{a})}
$$

To define quantum permutation relations we will use the
dynamical flow \eqref{6.3} as a {\it deforming flow\/}
considering the time variable as a quantum 
deformation parameter.

Denote $\lambda^\hb(A)=\rho(A)-\rho\big(\gamma^{-\hb}(A)\big)$, where
$\hb>0$, and set the permutation relations between
quantum generators:
\begin{equation}
[\CC,\BB]=\lambda^\hb(\aA),\quad 
\CC\aA=\gamma^\hb(\aA)\CC,\quad 
[\aA_\mu,\aA_\mu]=0,\quad 
\CC=\BB^*,\quad
\aA=\aA^*.
\label{6.5}
\end{equation}

Note that $\rho(\aA)-\CC\BB$ and $\varkappa_j(\aA)$
are the Casimir elements (the center elements) 
of the algebra generated by relations \eqref{6.5}. 
Also note that in the limit as $\hb\to0$, under the assumption
that $\frac{i}{\hb}[\,\cdot,\cdot\,]\to\{\,\cdot,\cdot\,\}$,  
relations \eqref{6.5} are transferred to \eqref{6.4} or (6.4a).

Let us introduce an associative quantum product $\star$ over
$\bR^{k+2}$ corresponding to algebra \eqref{6.5}:
\begin{equation}
[f(\stackrel{3}{\BB},\stackrel{2}{\aA},\stackrel{1}{\CC})]\cdot
[g(\stackrel{3}{\BB},\stackrel{2}{\aA},\stackrel{1}{\CC})]=
k(\stackrel{3}{\BB},\stackrel{2}{\aA},\stackrel{1}{\CC}),\quad
k=f\star g.
\label{6.6}
\end{equation}
Here $f,g$ are arbitrary polynomial. We use the normal
ordering of generators and pose indices $1,2,\dots$ indicating
their order from right to left as in [37]. 
The product $\star$ can be explicitly calculated using the
technique of ``regular representation'' (see, for example,
[27], Appendix 2). Namely,  
\begin{equation}
(f\star g)(B,A,C)=f(\stackrel{3}{L_B},\stackrel{2}{L_A},
\stackrel{1}{L_C})\,g(B,A,C),
\label{6.7}
\end{equation}
where $L_B=B\star$, $L_A=A\star$, $L_C=C\star$ are operators of
the left regular representation of permutation relations.
In the case of relations \eqref{6.5} it is rather easy to derive
(see [29]):
\begin{equation}
L_B=B,\quad L_A=\gamma^{\hb B\pa/\pa B}(A),\quad
L_C=C\gamma^{\hb\star}
+\Lambda^\hb (\stackrel{\vphantom{2}}{A},
\stackrel{2}{B}\lower1pt\hbox{$\stackrel{1}{{\pa}/{\pa B}}$} )\,
{\pa}/{\pa B}.
\label{6.8}
\end{equation}
Here the function $\Lambda^\hb$ is defined by 
$$
\Lambda^\hb(A,t)=\frac{\exp\{ t(\gamma^{\hb*}-I)\}-I }
{t(\gamma^{\hb*}-I)} \lambda^\hb(A),
$$
and $\gamma^{\hb*}$ denotes the shift operator by the variable
$A$, i.e., $(\gamma^{\hb*}f)(A)=f\big(\gamma^{\hb}(A)\big)$.

\begin{lemma} {}Formulas \eqref{6.7}, \eqref{6.8} determine the 
associative product $\star$ over $\bR^{k+2}$ satisfying 
\eqref{6.6}. 
This is the normal product with respect to the polarization 
$\cF^+=\{f(A,C)\}$ and the vacuum point $S_1=S_2=0$, 
$A=a$, where $a\in\bR^k$ is arbitrarily fixed.
\end{lemma}

Later on we denote the vacuum point 
$(0,0,a)\in\bR^{k+2}$ by the same letter~$a$,
as we did in Section~5.

\begin{example}
\rm
The Lie algebra $\su(1,1)$ is generated by the relations
$$
[\bs_1,\bs_2]=i\hb\bs_3,\qquad
[\bs_2,\bs_3]=-i\hb\bs_1,\qquad
[\bs_3,\bs_1]=-i\hb\bs_2.
$$
Denote $\CC=\bs_1+i\bs_2=\BB^*$, $\aA=\bs_3+\hb/2$. 
Then we obtain relations of the type \eqref{6.5}:
$$
\CC\BB=\BB\CC+2\hb\aA-\hb^2,\qquad \CC\aA=(\aA+\hb)\CC.
$$
The deforming flow $\gamma^t:\bR\to\bR$ is given by  
$\gamma^t(A)=A+t$, and 
$$
\lambda^\hb(A)=2\hb A-\hb^2,\qquad \rho(A)=A^2,\qquad 
\Lambda^\hb(A,t)=2\hb A+\hb^2 t-\hb^2.
$$
So, by formulas \eqref{6.8} we obtain
$$
L_B=B,\quad 
L_A=A+\hb\frac{\pa}{\pa B},\quad 
L_C=Ce^{\hb\pa/\pa A}+2\hb A\frac{\pa}{\pa B}
+\hb^2 B\frac{\pa^2}{\pa B^2}-\hb^2\frac{\pa}{\pa B}.
$$
Thus the normal $\star$-product corresponding to the 
polarization $\cF^+=\{f(A,C)\}$ of the enveloping 
of the Lie algebra $\su(1,1)$ is
determined by formula \eqref{6.7} with operators
$L_B,\,L_A,\,L_C$ given above.  
\end{example}

Now let us return to the general algebra \eqref{6.5}.
We define the quantum Wick product $*$ on symplectic 
leaves in $\bR^{k+2}$, the operation on the quantum
restriction to leaves, and describe the corresponding special
functions.  

The symplectic leaf $\cX_a\subset\bR^{k+2}$ passing through the
vacuum point (see Lemma~6.2) is determined by the following
values of constants in \eqref{6.2}: 
$K=\rho(a)$, $\const_j=\varkappa_j(a)$; that is,
\begin{equation}
\cX_a=\big\{S^2_1+S^2_2=\rho(A)-\rho(a),\,
\varkappa_j(A)=\varkappa_j(a)\big\}.
\label{6.9}
\end{equation}
Here we assume that $a\in\bR^k$ is not a point of local minimum
or maximum of the function $\rho$. Also we consider only the
connected component of the set determined by equations
\eqref{6.9}.
 
In order to introduce a complex structure on the leaf $\cX_a$
following the approach of Section~5, one has to fix a function
$Z=Z(A,C)\in\cF^+\setminus\cF^0$.
Let us choose a polynomial $g$ which is an integral of motion
for the flow $\gamma^t$ \eqref{6.3}, and $g(a)=\rho(a)$.
Assume that $\rho(A)>g(A)$ for all $A\ne a$ from a neighborhood
of the point~$a$.  
Now we split the polynomial $\rho(A)-g(A)$ in two polynomial
multipliers: 
\begin{equation}
\rho(A)-g(A)=D(A)\cdot E(A)
\label{6.10}
\end{equation}
in such a way that 
\begin{equation}
D(a)\ne0,\qquad E(a)=0.
\label{6.11}
\end{equation}
Then in the neighborhood of the point $A=a$
we define $Z(A,C)=C/\oD(A)$. 

One can parameterize points $A$ along the leaf \eqref{6.9} by the time 
on the trajectory: $A=\gamma^t(a)$. 
Then we have either the case $\rho\big(\gamma^t(a)\big)>\rho(a)$ 
for all $t>0$, or the case
\begin{equation}
\begin{aligned}
\rho\big(\gamma^t(a)\big)&>\rho(a)\qquad 
\text{for}\quad 0<t<t^*,\\
\rho\big(\gamma^{t^*}(a)\big)&=\rho(a).
\end{aligned}
\label{6.12}
\end{equation}
In the first case, we formally set $t^*=\infty$.

In the last case denote $a^*=\gamma^{t^*}(a)$. 
The point $S_1=S_2=0$, $A=a^*$ is ``polar'' with respect to the
vacuum point on the leaf $\cX_a$. 
In this case, in addition to \eqref{6.11}, we assume that the
polynomial~$D$ has zero at $a^*$:
\begin{equation}
D(a^*)=0.
\label{6.13}
\end{equation}
The equation 
\begin{equation}
z=C/\oD(A)
\label{6.14}
\end{equation}
determines the complex coordinate $z$ all over the leaf $\cX_a$
except the polar point. The polar point corresponds to $z=\infty$. 
In a neighborhood of the polar point one has to make the usual
change of variables taking the new complex coordinate 
$z'=1/z$.

Thus we introduce the global complex structure on the
symplectic leaf $\cX_a$ in both cases: 
$t^*=\infty$ and $t^*<\infty$. 

Note that in the case $t^*<\infty$ the leaf $\cX_a$ is compact 
and diffeomorphic to $\bS^2$. 
The value $t^*=t^*(a)$ depends on~$a$. We assume that~$a$ is
chosen in such a way that 
\begin{equation}
t^*(a)=(N+1)\hb,\qquad N\in\bZ_+.
\label{6.15}
\end{equation}

Now let us denote
$$
\cD(t)=D\big(\gamma^t(a)\big), \qquad
\cE(t)=E\big(\gamma^t(a)\big)
$$
and consider the following differential equation for the
function $k=k(r)$, $r\geq0$:
\begin{equation}
\overline{\cE}\Big(\hb r\frac{d}{dr}\Big)k
=r\cD\Big(\hb r\frac{d}{dr}+\hb\Big)k, \qquad k(0)=1.
\label{6.16}
\end{equation}
The solution is given by the series
\begin{equation}
k(r)=1+\sum^{N}_{n=1}\frac{r^n}
{\cH(\hb)\cH(2\hb)\dots\cH(n\hb)},\qquad
\cH\,\od\,\frac{\overline{\cE}}{\cD}
\label{6.17}
\end{equation}
Here $N=\infty$ if $t^*=\infty$, and we assume that the
multipliers $D$, $E$ are taken in such a way that the series
\eqref{6.17} has a convergency domain $\{|z|<\infty\}$.

Since $\cD$ and $\cE$ are polynomial, the function $k$
\eqref{6.17} is of hypergeometric type. Also $k>0$ everywhere.
So one can define the K\"ahler form on the leaf~$\cX_a$ by the
formula: 
\begin{equation}
\omega\,\od\, i\hb\opa\pa\ln k(|z|^2).
\label{6.18}
\end{equation}
We call $\omega$ the {\it hypergeometric K\"ahler form\/} on the
surface of revolution $\cX_a$ \eqref{6.9} with the complex
structure \eqref{6.14}. 

\begin{lemma}
{\rm(a)} The form \eqref{6.17} is globally defined on $\cX_a$,
and condition \eqref{6.15} is equivalent to \eqref{2.2}{\rm:}
$\frac1{2\pi\hb}\int_{\cX_a}\omega=N$.

{\rm(b)} The function $\cK(\oz|z)=k(|z|^2)$ is the reproducing
kernel for the space $\cL_a$ of antiholomorphic distributions
with the Hilbert norm 
$$
\|u\|=\Big(\sum_{n\geq0}\cH(\hb)\dots\cH(n\hb)|u_n|^2\Big)^{1/2},
\qquad \text{where}\quad u(\oz)=\sum_{n\geq0}\oz^nu_n.
$$
This is the same norm and the same reproducing kernel as could
be obtained by general formulas \eqref{5.6}, \eqref{5.5}{\rm;}
the Darboux coordinate $P$ dual to $Z=C/\oD(A)$ 
{\rm(}see \eqref{5.1},\eqref{5.2}{\rm)} 
in this case is given by the formula $P=Bt(A)/\overline{E}(A)$,
where $t=t(A)$ is time along trajectories of the flow
\eqref{6.3} with the initial data $t(a)=0$.

{\rm(c)} The asymptotics of $\cK$ as $\hb\to0$ is the
following{\rm:} 
$$
\cK(\oz|z)=\const\sqrt{g_0(|z|^2)}\, e^{F_0(|z|^2)/\hb}(1+O(\hb)),
\qquad g_0(r)=\big(rF'_0(r)\big)',
$$
where $F_0$ is the solution of the problem
$$
\cH\big(rF'_0(r)\big)=0,\qquad F_0(0)=0,
$$
and $\const=\lambda(a)^{1/2}/|D(a)|$. The hypergeometric form
\eqref{6.18} has the asymptotics{\rm:} 
$$
\omega=\omega_0+O(\hb),\qquad
\omega_0\,\od\,ig_0(|z|^2)\,d\oz\wedge dz,
$$
where $\omega_0$ is the classical symplectic form on the
symplectic leaf $\cX_a$ generated from $\bR^{k+2}$ by the
Poisson structure~\eqref{6.4}. 

{\rm(d)} The Hilbert norm in $\cL_a$ can be represented in
the integral form \eqref{5.7} via the reproducing measure  
\begin{equation}
dm=k(|z|^2)\ell(|z|^2)\,d\oz dz
\label{6.19}
\end{equation}
if there exists a positive solution of the hypergeometric equation{\rm:}
\begin{equation}
\overline{\cE}\Big(-\hb r\frac{d}{dr}\Big)\ell
=r\cD\Big(-\hb r\frac{d}{dr}-\hb\Big)\ell,\qquad
\frac1{\hb}\int^\infty_0\ell(r)\,dr=1.
\label{6.20}
\end{equation}
The asymptotics is the following{\rm:} $dm=dm^{\omega_0}+O(\hb)$.
\end{lemma}

These statements were proved in [29]. 
Pay attention to equation \eqref{6.20}, which differs from
\eqref{6.16} by changing the sign near the parameter~$\hb$ only.
But the properties of the function~$\ell$ in \eqref{6.20} 
are certainly not the same as of the function~$k$ in
\eqref{6.16}; 
for instance, $\ell$ can have a weak singularity at $r=0$, 
and it is decreasing as $r\to\infty$. A series of examples was
considered in [29].

\hbox to 1cm{}

Now let us come back to the general Theorem~5.1. 
In the given example the irreducible representation \eqref{5.9} 
of the algebra \eqref{6.5} can be realized in the space
$\cL_a=\cL(\cX_a)$ by the following operators{\rm:}
\begin{equation}
\wh{A}_a=\gamma^{\hb\oz\opa+\hb}(a),\qquad 
\wh{B}_a=\cD(\hb\oz\opa)\cdot\oz,\qquad 
\wh{C}_a=\frac1{\oz}\cdot\cE(\hb\oz\opa).
\label{6.21}
\end{equation}
Symbols $A_a,\,B_a,\,C_a$ are the quantum restriction of the
coordinate functions $A,\,B,\,C$ to the leaf $\cX_a$.
They can be calculated, say, by \eqref{3.9}{\rm:} 
\begin{align*}
A_a&=\frac1{\cK}\wh{A}_a(\cK)=\gamma^{\hb rd/dr+rF'(r)+\hb}(a)1(r),\\
C_a&=\frac1{\cK}\wh{C}_a(\cK)=\frac1{\oz}
\cE\big(\hb rd/dr+rF'(r)\big)1(r)=\oB_a,
\end{align*}
where $r=|z|^2$ and $F$ is the quantum K\"ahler potential 
$F(r)=\hb\ln k(r)$.

\begin{theorem}
Let $\star$ be the normal product over $\bR^{k+2}$ defined by
\eqref{6.7}, \eqref{6.8}, and let $*$ be the Wick product over
the surface of revolution $\cX_a\subset\bR^{k+2}$ generated by
the hypergeometric K\"ahler form \eqref{6.18}. 
Then for any polynomial $f=f(B,A,C)$ on $\bR^{k+2}$ the quantum 
restriction to the surface~$\cX_a$ is given by  
$$
f\big|\zs{\wh{\cX}_a}\equiv f_a=\frac1{\cK}f
\Big(\stackrel{3}{\wh{B}}_a,\stackrel{2}{\wh{A}}_a,
\stackrel{1}{\wh{C}}_a\Big)(\cK),
$$
where the differential operators
$\wh{B}_a,\,\wh{A}_a,\,\wh{C}_a$ are defined 
in \eqref{6.21} and $\cK=k(\oz z)$ is the hypergeometric
function \eqref{6.17}.  
The conditions \eqref{3}, \eqref{4} are satisfied.
\end{theorem}

Note that the key role in this construction is played by the 
hypergeometric reproducing kernel $\cK$ which can be also 
considered as the coherent state in $\cL(\cX_a)$ corresponding
to the irreducible representation \eqref{6.21} of the algebra
\eqref{6.5}:
\begin{equation} 
\cK(\oz|z)=\exp\Big\{\frac{z\wh{P}_a}\hb\Big\}1(\oz)
=\bigg(I+\sum_{n\geq1}\frac{z^n}
{\overline{\cE}(n\hb)\dots\overline{\cE}(\hb)}\wh{B}^{\,n}_a \bigg)1(\oz).
\label{6.22}
\end{equation}
Here $1=1(\oz)$ is the ``vacuum'' element in $\cL(\cX_a)$, 
the operator $\wh{B}_a$ is given in \eqref{6.21}, 
and 
$$
\wh{P}_a=\wh{B}_a\,\overline{E}^{-1}(\wh{A}_a)t(\wh{A}_a)
=\hb\oz\cH(\hb\oz\opa+\hb)^{-1}(\oz\opa+1).
$$
The operator acting to the vacuum element in \eqref{6.22} 
is a hypergeometric function in
the creation operator $\wh{B}_a$, or just the exponential
function in the creation operator~$\wh{P}_a$.
This second creation operator $\wh{P}_a$ is not very convenient
since, in general, it is pseudodifferential. 
So, staying in the class of representations of
algebra~\eqref{6.5} by differential operators, 
we have to deal with hypergeometric functions instead of 
the standard exponential function.

\setcounter{example}{0}

\begin{example} \rm (continuation).
In the case of the Lie algebra $\su(1,1)$, 
the symplectic leaves \eqref{6.9} are hyperboloids: 
$\cX_a=\{A^2-BC=a^2,A\geq a\}$, where $a=\const>0$.
Let us fix $g(A)=\const=\rho(a)=a^2$. 
There are two possible choices of the multipliers
$D,\,E$ \eqref{6.10} which factorize the difference 
$\rho(A)-\rho(a)=(A-a)(A+a)$.

{\bf Variant I}. {}First, we can take $E(A)=A-a$ and $D(A)=A+a$.  
Then $\cE(t)=t$, $\cD(t)=t+2a$, and equations \eqref{6.16},
\eqref{6.20} read 
\begin{equation}
\hb(1-r)dk/dr=(2a+\hb)k,\qquad
\hb(1-r)d\ell/dr=(\hb-2a)\ell.
\label{6.23}
\end{equation}
Since $k(0)=1$, then $k(r)=(1-r)^{-(2a+\hb)/\hb}$. 
This function is singular at $r=1$, what means 
that the convergency radius of the series \eqref{6.17} in this 
variant is equal to~$1$ (but not to $\infty$).
Hence the normalization condition for the function $\ell$ is 
$\frac1\hb\int^1_0\ell(r)\,dr=1$. 
Solving the second equation \eqref{6.23}, we obtain 
$\ell(r)=2a(1-r)^{(2a-\hb)/\hb}$. Thus the quantum symplectic
form and the reproducing measure in this variant are given by 
$$
\omega=i(2a+\hb)\frac{d\oz\wedge dz}{(1-|z|^2)^2},\qquad
dm=2a\frac{d\oz dz}{(1-|z|^2)^2}.
$$
They are transported to the hyperboloid $\cX_a$ by means of the
$\su(1,1)$-invariant complex structure \eqref{6.14}: 
$$
z=\frac{C}{A+a}\quad\text{or}\quad
A=\frac{a(1+|z|^2)}{1-|z|^2},\quad
C=\frac{2az}{1-|z|^2}.
$$
Of course, the quantum geometrical data $\omega$, $dm$ are
$\hb$-deformations of the classical data $\omega_0$,
$dm^{\omega_0}$, where $\omega_0$ is the $\su(1,1)$-invariant
K\"ahlerian form on the hyperboloid: 
$\omega_0=2ia(1-|z|^2)^{-2}d\oz\wedge dz
=\frac{i}2A^{-1}dB\wedge dC\big|\zs{\cX_a}$.
In this variant, all the operators \eqref{6.21} are of the first
order and generate the standard irreducible representation of
the Lie algebra $\su(1,1)$.

{\bf Variant II}. Let us make another choice: 
$E(A)=A^2-a^2$ and $D(A)=1$. Then $\cE(t)=t(2a+t)$ and, instead
of \eqref{6.23}, we have the following equations  
for the functions~$k$ and~$\ell$:
\begin{align*}
&\hb^2 rd^2k/dr^2+(2a\hb+\hb^2)dk/dr-k=0,\\
&\hb^2 rd^2\ell/dr^2-(2a\hb-\hb^2)d\ell/dr-\ell=0.
\end{align*}
The normalization conditions are the same as in \eqref{6.16},
\eqref{6.20}; so the solution is
\begin{equation}
k(r)=\wt{I}_{2a/\hb}\Big(\frac{2}{\hb}\sqrt{r}\Big),\qquad 
\ell(r)=\wt{M}_{2a/\hb}\Big(\frac{2}{\hb}\sqrt{r}\Big).
\label{6.24}
\end{equation}
Here $\wt{I}_\nu$ and $\wt{M}_\nu$ are modified Bessel and MacDonald
functions: 
\begin{align*}
\wt{I}_\nu(y)\,\od\,\sum_{n\geq0}\Big(\frac{y}{2}\Big)^{2n}
\frac{\Gamma(\nu+1)}{n!\Gamma(\nu+n+1)},
\\
\wt{M}_\nu(y)\,\od\,\frac{(y/2)^\nu}{\hb\Gamma(\nu+1)}
\int^{\infty}_{-\infty}\exp\{-y\cosh t-\nu t\}\,dt.
\end{align*}
Thus the quantum symplectic form and the reproducing measure in 
this variant are given by 
\begin{equation}
\omega=i\hb\opa\pa\ln\wt{I}_{2a/\hb}(2|z|/\hb),
\qquad
dm=\big(\wt{I}_{2a/\hb}\wt{M}_{2a/\hb}\big)\big(2|z|/\hb\big)\,
d\oz dz.
\label{6.25}
\end{equation}
They are transported to the hyperboloid $\cX_a$ by means of the 
complex structure \eqref{6.14}:
\begin{equation}
z=C,\qquad A=(a^2+|z|^2)^{1/2}.
\label{6.26}
\end{equation}
The operators \eqref{6.21} in this variant are the following:
\begin{equation}
\wh{A}_a=a+\hb+\hb\oz\opa,\qquad
\wh{B}_a=\oz,\qquad
\wh{C}_a=\hb\opa\cdot(2a+\hb\oz\opa).
\label{6.27}
\end{equation}
Since the complex structure \eqref{6.26} is not
$\su(1,1)$-invariant, these operators are not all of the first
order; they generate the irreducible representation of
$\su(1,1)$ given in [1].  
The corresponding ``Bessel'' geometrical data \eqref{6.25} on
the hyperboloid are quantum $\hb$-deformations of the classical 
data $\omega_0$, $dm^{\omega_0}$, where 
$\omega_0=\frac{i}2(a^2+|z|^2)^{-1/2}\,d\oz\wedge dz
=\frac{i}2A^{-1} dB\wedge dC\big|\zs{\cX_a}$.
\end{example}

\begin{example} {\bf Quadratic algebra of the Zeeman effect}.
\rm  
The Hydrogen atom in a homogeneous magnetic field is described
by the Hamiltonian:
$\bH=\big(\wh{p}-\cA(q)\big)^2-|q|^{-1}$, 
where $\wh{p}=-i\hb\pa/\pa q$, $q\in\bR^3$, 
and the magnetic potential $\cA$ has the following components 
$\cA_1=-\frac12\ve q_2$, $\cA_2=\frac12\ve q_1$, $\cA_3=0$. 
Here we assume that the magnetic field is directed along the
third coordinate axis, and the values of~$\ve$ characterize  
the strength of the field. 
Applying the quantum averaging procedure, we can transform  
this Hamiltonian with arbitrary accuracy $O(\ve^{n+2})$ to the
following form 
\begin{equation}
\bH\sim\bH_0+\ve \bM_3+\ve^2 
f^{(n)}(\bbS_0,\bbS_1,\bbS_2,\bbS_3;\bH_0,\bM_3)
+O(\ve^{n+2}).
\label{6.28}
\end{equation}
Here $\bH_0=|q|(\wh{p}^2+\frac14)$ describes the Hydrogen atom
itself, $\bM_3$ is the third component of the angular momentum
$\bM=q\times\wh{p}$, and the operators $\bbS_j$ 
generate the {\it algebra of joint symmetries} of $\bH_0$ and $\bM_3$,
that is, $[\bH_0,\bbS_j]=[\bM_3,\bbS_j]=0$, $j=0,\dots,4$.
More precisely, formulas for $\bbS_j$ are
$$
\bbS_0=\bL_3-\bbR_3,\quad
\bbS_1=\bL_1\bbR_2-\bL_2\bbR_1,\quad
\bbS_2=\bL_1\bbR_1+\bL_2\bbR_2,\quad
\bbS_3=\bL_3\bbR_3+\bL^2,
$$
where $\bL=(\bM+\bN)/2$, $\bbR=(\bM-\bN)/2$, and 
$\bN=q(\wh{p}^2+\frac14)-\wh{p}\times\wh{\bM}+\wh{\bM}\times\wh{p}$.

On the right-hand side of \eqref{6.28} we assume that 
generators $\bs_j$ are ordered in some way. 
The function $f^{(n)}$ is a polynomial of
degree~$n$ in the variables~$S_j$; 
details and explicit formulas for $f^{(2)}$ see in [28]. 

The most important fact is that the generators $\bbS_j$
satisfy the following quadratic permutation relations:
\begin{align}
[\bbS_1,\bbS_2]&=\frac{i\hb}{2}(\bbS_0\bbS_3+\bbS_3\bbS_0),&\qquad 
[\bbS_0,\bbS_1]&=2i\hb\bbS_2,\nonumber\\
[\bbS_2,\bbS_3]&=-\frac{i\hb}{2}(\bbS_0\bbS_1+\bbS_1\bbS_0),&\qquad 
[\bbS_0,\bbS_2]&=-2i\hb\bbS_1,
\label{6.29}\\
[\bbS_3,\bbS_1]&=-\frac{i\hb}{2}(\bbS_0\bbS_2+\bbS_2\bbS_0),&\qquad 
[\bbS_0,\bbS_3]&=0.
\nonumber
\end{align}

Let us introduce $\CC=\bbS_1+i\bbS_2$, $\BB=\bbS_1-i\bbS_2$,  
and also denote $\aA_1=\bbS_0-\hb$,
$\aA_2=\bbS_3+\frac\hb2\bbS_0$.  
Then the relations \eqref{6.29} can be written in the form
\eqref{6.5}, where $\rho(A)\,\od\,A^2_2$ and the dynamical flow
\eqref{6.3} in $\bR^2$ 
is generated by the vector field 
$v=2\pa/\pa A_1-A_1\pa/\pa A_2$.  
The integral of motion of this field is
$\varkappa(A)=A^2_1+4A_2$. Thus the symplectic leaf \eqref{6.9}
is defined as follows: 
\begin{equation}
\cX_a=\{A^2_2-BC=a^2_2,\, A^2_1+4A_2=a^2_1+4a_2\}.
\label{6.30}
\end{equation}
The representation of the algebra \eqref{6.29} related to the
Zeeman effect selects the values of the parameters $a_1$, $a_2$:
namely, $a_1<0$, $a_2>0$.
Then the leaf $\cX_a$ \eqref{6.30} is topologically a sphere.

Evaluating 
$\rho\big(\gamma^t(a)\big)-\rho(a)=t(t-|a_1|)(t-t_+)(t-t_-)$,
where $t_{\pm}=\frac12|a_1|\pm\frac12(a^2_1+8a_2)^{1/2}$, 
we conclude that the inequality \eqref{6.12} is satisfied on the
interval $0<t<t^*=|a_1|$. Taking into account \eqref{6.15}, we
obtain the quantization condition 
$$
a_1=-(N+1)\hb,\qquad N\in\bZ_+.
$$
The polar point $a^*=\gamma^{t^*}(a)$ is the following: 
$a^*=(-a_1,a_2)$.

Now let us choose $g(A)=\frac14t^2_+(A^2_1+4A_2-t^2_+)$. 
Then the polynomial $\rho(A)-g(A)$ can be factorized,
as in \eqref{6.10}, taking
$$
D(A)=A_2+\frac12t_+A_1-a_2+\frac12t_+a_1,\qquad
E(A)=A_2-\frac12t_+A_1-a_2+\frac12t_+a_1.
$$
So, we have
\begin{equation}
\cD(t)=D\big(\gamma^t(a)\big)=(a_1+t)(t_+-t),\quad
\cE(t)=E\big(\gamma^t(a)\big)=t(t_--t).
\label{6.31}
\end{equation}
Thus the hypergeometric equations \eqref{6.16}, \eqref{6.20}
for~$k$ and for~$\ell$ are the following:
\begin{align*}
&r(1-r)\frac{d^2k}{dr^2}+\big(\alpha_+-(\alpha_--N+1)r\big)
\frac{dk}{dr}+N\alpha_-k=0,\\
&r(1-r)\frac{d^2\ell}{dr^2}+\big(\beta_--(\beta_++N+3)r\big)
\frac{d\ell}{dr}-(N+2)\beta_+\ell=0,
\end{align*}
where $\alpha_{\pm}=1-t_{\mp}/\hb$, $\beta_{\pm}=1+t_{\pm}/\hb$.
The solution of these equations could be expressed
via the Gauss hypergeometric functions 
${}_2F_1(-N,\alpha_-,\alpha_+;r)$ and 
${}_2F_1(N+2,\beta_+,\beta_-;r)$. But we prefer
to demonstrate explicit formulas:
\begin{align}
k(r)&=\sum^N_{n=0}\frac{N!}{n!(N-n)!}
\frac{(t_+-\hb)(t_+-2\hb)\dots(t_+-n\hb)}
{(|t_-|+\hb)(|t_-|+2\hb)\dots(|t_-|+n\hb)}r^n,
\label{6.32}\\
\ell(r)&=
\frac{\hb(N+1)\Gamma(1+(t_++|t_-|)/\hb)}
{\Gamma(t_+/\hb)\Gamma(1+|t_-|/\hb)}
\int^\infty_0
\frac{\lambda^{t_+/\hb}\,d\lambda}
{(1+\lambda r)^{N+2}(1+\lambda)^{1+(t_++|t_-|)/\hb}}.
\nonumber
\end{align}
Here $\Gamma$ is the gamma-function, 
and in the last integral representation the Euler formula for
hypergeometric functions was used.  

We see that in this example the reproducing kernel 
$\cK(\oz|z)=k(\oz z)$ is a {\it Jacobi polynomial\/}
\eqref{6.32} of degree~$N$.  
{}Functions \eqref{6.32} generate the quantum K\"ahlerian form
and the reproducing measure on the ``sphere'' $\cX_a$ \eqref{6.30} by
the general formulas \eqref{6.18}, \eqref{6.19}. The complex
structure on $\cX_a$ is determined by \eqref{6.14}:
$z=C/(A_2+\frac12t_+A_1-a_2+\frac12t_+a_1)$.
The corresponding irreducible representation of the algebra
\eqref{6.29} in the space $\cL(\cX_a)$ is given by the
operators: 
$$
\wh{S}_1=\frac12(\wh{C}+\wh{B}),\quad
\wh{S}_2=\frac1{2i}(\wh{C}-\wh{B}),\quad
\wh{S}_3=\wh{A}_2-\frac\hb2\wh{A}_1-\frac{\hb^2}{2},\quad
\wh{S}_0=\wh{A}_1+\hb,
$$
where $\wh{C},\,\wh{B},\,\wh{A}$ are defined via \eqref{6.21},
\eqref{6.31}: 
\begin{align*}
\wh{A}_1&=a_1+2\hb+2\hb\oz\opa,&\quad
\wh{A}_2&=a_2-(a_1\hb+a_1\hb\oz\opa)-(\hb\oz\opa+\hb)^2,\\
\wh{B}&=(a_1+\hb\oz\opa)(t_+-\hb\oz\opa)\cdot\oz,&\quad
\wh{C}&=\hb\opa\cdot(t_--\hb\oz\opa).
\end{align*}

This is the realization of relations \eqref{6.29} by the second
order differential operators. The equivalent representations can
be obtained by three other possible variants of the factorization
\eqref{6.10} (or choices of complex structures):
\begin{align*}
\cD&=(a_1+t)(t_--t)&\,\cD&=(a_1+t)(t_+-t)(t_--t)&\,\cD&=a_1+t\\
\cE&=t(t_+-t),      &\,\cE&=t,                    &\,\cE&=t(t_+-t)(t_--t).
\end{align*}
So, together with the variant considered above, 
totally we have four different algebraic complex structures on
the symplectic leaf $\cX_a\subset \bR^4$. 
They generate four quantum K\"ahler structures 
and four irreducible equivalent representations of the algebra
\eqref{6.29}  by differential operators (up to order six) [29].

Other examples and a list of references around algebras of the type
\eqref{6.5} see in [29].

\end{example}

\setcounter{equation}{0}

\section{Quantum cylinder and theta-functions}

{}Finally, we consider surfaces of revolution \eqref{6.2}
\begin{equation}
\cX=\{\rho(A)-BC=K,\,\varkappa_j(A)=\const_j \,(j=1,\dots,k-1)\}
\label{7.1}
\end{equation}
under the condition $K<\min\rho$. In this case, the values 
$BC=|C|^2$ are strictly positive on $\cX$, and so, no vacuum
point exists. Instead of that there are noncontractible
$1$-cycles: sections of $\cX$ by the planes $A=\const$.
Thus the surface $\cX$ topologically is the cylinder or the
torus. Actually, the axis of revolution for the surface $\cX$ is
determined by the trajectory $\{\gamma^t(a_0)\}$,
starting from a point $a_0$, where $\varkappa_j(a_0)=\const_j$ 
$(j=1,\dots,k-1)$. 
If this trajectory is periodic, then $\cX$ is homeomorphic to
the torus $\bT^2$;
if the trajectory is not periodic, then $\cX$ is homeomorphic to
the cylinder $\bS\times\bR$.
Here we consider only the cylinder case.

As above, let us take some polynomial $g(A)$ such that
$g(\gamma^t(a_0))=g(a_0)=K$ for any $t$, and $\rho(A)>g(A)$ 
in a neighborhood of $\cX$. Moreover, we assume that there is a
(complex) polynomial $M(A)$ such that 
\begin{equation}
\rho(A)-g(A)=|M(A)|^2.
\label{7.2}
\end{equation}
Then one can choose the following multipliers in the 
factorization \eqref{6.10}: 
$D(A)=M(A)e^{-t(A)}$, $E(A)=\overline{M}(A)e^{t(A)}$, where
$t(A)$ is the time along trajectories of the flow \eqref{6.3}.

Since instead of the vacuum point on $\cX$ there is a
noncontractible circle, it is natural to consider a ring in the
complex plane as a coordinate chart. 
That is why in all formulas \eqref{6.14}, \eqref{6.18}, 
\eqref{6.21}, \eqref{6.22} 
we now replace the complex variable $z$ by $e^z$, 
and the real variable $r=z\oz$ we replace by $r=z+\oz$.
Then the operator $\hb r\,d/dr$ played the role of 
the ``quantum time'' in \eqref{6.16}, \eqref{6.20} must be
replaced by the operator $\hb\,d/dr$.

In order to obtain an analog of the equation \eqref{6.16} in the
cylinder case, let us first remark that \eqref{6.16} is equivalent to
the equations  
\begin{equation}
\wh{A}\cK=\overline{\wh{A}}\cK,\qquad
\wh{B}\cK=\overline{\wh{C}}\cK.
\label{7.3}
\end{equation}
Here the bar $\overline{\vphantom{a}\dots}$ denotes the complex
conjugation of an operator;  
after this conjugation the operator acts by~$z$ (not by~$\oz$). 
The operators $\wh{A},\,\wh{B},\,\wh{C}$ in \eqref{7.3} are
generators of the irreducible representation of algebra
\eqref{6.5} in the space of antiholomorphic sections. 
Now we have to use not the operators of type \eqref{6.21} but
the following ones: 
$$
\wh{A}=\gamma^{\hb\opa+\hb}(a_0),\qquad
\wh{B}=\cD(\hb\opa)\cdot e^{\oz},\qquad
\wh{C}=e^{-\oz}\cdot\cE(\hb\opa).
$$
Thus equations \eqref{7.3} imply
$$
\cK(\oz|z)=k(\oz+z),\qquad 
M(\hb\opa)e^{-\hb\opa}e^{\oz}k(\oz+z)
=e^{-z}M(\hb\pa)e^{\hb\pa}k(\oz+z).
$$
Since $M\ne0$, we obtain the equations which do not depend on
the polynomial~$M$ at all:
\begin{equation}
k(r+2\hb)=e^rk(r),\qquad k(r+2\pi i)=k(r).
\label{7.4}
\end{equation}
The last condition follows form the $2\pi i$-periodicity 
of $\cK(\oz|z)$ by both~$\oz$ and~$z$.

The equations for the function~$\ell$ (the density of the
reproducing measure) in general are the following:
$$
\wh{A}\,{\vphantom{\Big|}}'\ell
=\overline{\wh{A}\,}\,{\vphantom{\Big|}}'\ell,
\qquad
\wh{B}\,{\vphantom{\Big|}}'\ell
=\overline{\wh{C}}\,{\vphantom{\Big|}}'\ell,
$$
where the prime $'$ denotes the transposition with respect to the
standard pairing of functions:
$\langle\varphi,\psi\rangle=\int\varphi\psi\,d\oz dz$.
{}From these equations, we obtain $\ell=\ell(r)$, 
where $r=\oz+z$, and 
\begin{equation}
\ell(r+2\hb)=e^{-r}\ell(r),\qquad \frac{1}{\hb}
\int^{\infty}_{-\infty}\ell(r)\,dr=1.
\label{7.5}
\end{equation}
The solution of the last problem is just the Gaussian function
\begin{equation}
\ell(r)=(\hb/4\pi)^{1/2}e^{-(r-\hb)^2/4\hb}.
\label{7.6}
\end{equation}
The solution of equations \eqref{7.4} is given by the series
\begin{equation}
k(r)=\sum_{n\in\bZ}e^{-n^2\hb+n(r-\hb)}=\theta(r-\hb,e^{-\hb}).
\label{7.7}
\end{equation}
Here we denote by $\theta$ the following theta-function:
\begin{equation}
\theta(\alpha,q)\,\od\,\sum_{n\in\bZ}q^{n^2}e^{n\alpha},\qquad 
q<1.
\label{7.8}
\end{equation}

\begin{theorem}
Assume that the trajectory $\{\gamma^t(a_0)\}$ is not periodic and
thus the surface $\cX$ \eqref{7.1} is homeomorphic to the
cylinder. Then  

{\rm(a)} In the space $\cL$ of antiholomorphic $2\pi i$-periodic
functions with the Hilbert norm
$$
\|u\|=\bigg(\frac1{4\pi\sqrt{\pi \hb}}\int_{0<\Im z<2\pi}
|u(\oz)|^2e^{-(\oz+z-\hb)^2/4\hb}\,d\oz dz\bigg)^{1/2},
$$
the irreducible representation of algebra \eqref{6.5} acts
by the following operators{\rm:}
\begin{equation}
\wh{A}=\gamma^{\hb\opa+\hb}(a_0),\quad
\wh{B}=e^{\oz-\hb\opa-\hb/2}M(\wh{A}),\quad
\wh{C}=\overline{M}(\wh{A})e^{\hb\opa-\oz+\hb/2},
\label{7.9}
\end{equation}
where the polynomial $M$ is taken from \eqref{7.2}.
This representation, up to equivalence, does not depend on the
choice of the point $a_0$, but is parametrized by values of the
Casimir functions $K$, $\varkappa_j$ in \eqref{7.1}.

{\rm(b)} The space $\cL$ is identified with the space $\cL(\cX)$
of antiholomorphic sections over the symplectic leaf \eqref{7.1} 
supplied with the complex structure 
$e^z=Ce^{t(A)}/\overline{M(A)}$ and with the K\"ahlerian form
\begin{equation}
\omega=i\hb\opa\pa\ln\cK,\qquad
\cK(\oz|z)=\theta(\oz+z-\hb,e^{-\hb}),
\label{7.10}
\end{equation}
where the theta-function $\theta$ is defined by \eqref{7.8}.

{\rm(c)} The reproducing measure on $\cX$ corresponding to the 
K\"ahlerian form \eqref{7.10} is given by 
\begin{equation}
dm=\frac12\theta\Big(\frac{i\pi}{\hb}(\oz+z-\hb),\,
e^{-\pi^2/\hb}\Big)\,d\oz dz.
\label{7.11}
\end{equation}
\end{theorem}

The last statement about the measure 
$dm=k(\oz+z)\ell(\oz+z)\,d\oz dz$ is derived from \eqref{7.6}, 
\eqref{7.7} and from the Jacobi transform of theta-functions
which reads in our case  
\begin{align}
\cK(\oz|z)&\equiv k(\oz+z)
=(\pi/\hb)^{1/2}e^{(\oz+z-\hb)^2/4\hb}
\theta\Big(\frac{i\pi(\oz+z-\hb)}{\hb},\,e^{-\pi^2/\hb}\Big)
\nonumber\\
&=(\pi/\hb)^{1/2}e^{(\oz+z-\hb)^2/4\hb}
\big(1+O(e^{-\pi^2/\hb})\big).
\label{7.12}
\end{align}

Now let us note that after the substitution of \eqref{7.12}
into the right-hand side of \eqref{7.10},  we obtain
\begin{align}
\omega&=\omega_0+i\hb\opa\pa\ln
\theta\Big(\frac{i\pi(\oz+z-\hb)}{\hb},e^{-\pi^2/\hb}\Big)
\nonumber\\
&=\omega_0+O(e^{-\pi^2/\hb})\qquad \text{as}\quad \hb\to0.
\label{7.13}
\end{align}
Here $\omega_0$ denotes the classical symplectic form on the
leaf $\cX$ \eqref{7.1}: 
$$
\omega_0=\frac{i}{2}d\oz\wedge dz.
$$
And also from \eqref{7.11} we obtain the asymptotics 
\begin{equation}
dm=dm^{\omega_0}+O(e^{-\pi^2/\hb}),\qquad\text{where}\quad
dm^{\omega_0}=\frac12d\oz dz.
\label{7.14}
\end{equation}

Let us stress that 
the exponentially small remainders $O(e^{-\pi^2/\hb})$
in \eqref{7.13}, \eqref{7.14}
are not zero and  {\it explicitly given  by the theta-series}. 
This means that the corresponding Wick product over the leaf
$\cX$ has the asymptotic expansion 
\begin{equation}
\psi*\chi\simeq \sum^{\infty}_{r=0}\frac{\hb^r}{r!}
\pa^r\! \psi\,\opa^{\,r}\!\chi + O(e^{-\pi^2/\hb}).
\label{7.15}
\end{equation}
The formal $\hb$-series on the right-hand side of \eqref{7.15}
is the well-known expansion of the Wick product over $\bR^2$
corresponding to the symplectic form $\omega_0$ and the measure
$dm^{\omega_0}$. 
This series does not know anything about the periodic condition 
on functions transforming the plane $\bR^2$ to the cylinder 
$\bR\times \bS$ (in our case: $2\pi i$-periodicity by the $z$
variable). 
Only the exponentially small remainder in \eqref{7.15} ``knows''
about this additional periodic condition, and so, about the
topology of $\cX$. This remainder is invisible within the frames
of  formal  constructions of $*$-products by power series
in $\hb$. 
In order to take into account the nontrivial $1$-cycle on $\cX$,
we have to exploit formulas for the Wick $*$-product including
the exponentially small quantities $\exp\{-j\pi^2/\hb\}$
for all $j=1,2,\dots$.
We consider this fact as a {\it tunneling} along the
noncontractible cycle on the quantum cylinder. 

The same results can be obtained in the case where $\cX$ is
homeomorphic to the torus $\bT^2=\bS\times\bS$ [30],
and obviously, this is a general property of the quantum
geometry over not simply connected K\"ahler manifolds.

\begin{example}
\rm
The simplest cylindric symplectic leaves are one-sheet
hyperboloids in $\bR^3=\su(1,1)^*$.
So, let us return to the Lie algebra $\fg=\su(1,1)$ already
considered in Example~6.1. 
Now we take the following family of symplectic leaves:
\begin{equation}
\cX=\{BC-A^2=\lambda^2\},\qquad \lambda>0.
\label{7.16}
\end{equation}
In this case $K=-\lambda^2$ and $\rho(A)=A^2$.
We choose $g(A)\equiv -\lambda^2$, $a_0=0$ and so 
$\rho(A)-g(A)=A^2+\lambda^2=(A+i\lambda)(A-i\lambda)$. 
Thus one can take $M(A)=A-i\lambda$ in \eqref{7.2} and
introduce the complex structure to the hyperboloid \eqref{7.16}:
\begin{equation}
e^z=C\frac{e^{A}}{A+i\lambda}\quad\text{or}\quad
A=\RR z,\quad 
C=(\RR z+i\lambda)e^{i\Im z}.
\label{7.17}
\end{equation}
The K\"ahlerian form and the reproducing measure on the
hyperboloid \eqref{7.16} are defined by the general formulas
\eqref{7.10},  \eqref{7.11}. 
The irreducible Hermitian representation of the Lie algebra
$\su(1,1)$ (the ``prime series'') is given by formulas
\eqref{7.9} which read: 
\begin{align}
\wh{A}=\hb+\hb\opa,\qquad
\wh{B}&=e^{\oz-\hb\opa-\hb/2}(\hb-i\lambda+\hb\opa),
\label{7.18}\\
\wh{C}&=(\hb+i\lambda+\hb\opa)e^{\hb\opa-\oz+\hb/2}.\nonumber
\end{align}
In our case, the representation contains pseudodifferential
operators since it corresponds to the complex structure
\eqref{7.17} which is not $\su(1,1)$-invariant. 
The Casimir element in this representation takes the value
$\wh{A}^2-\wh{C}\wh{B}=-\lambda^2\cdot I$. 
Usually, the prime series is realized by differential operators
on the circle what corresponds to the choice of the real
polarization over the hyperboloid $\cX$ \eqref{7.16}.
\end{example}

Note that the use of a real polarization over not simply
connected symplectic leaves generates associative algebras of functions 
which are not central (that is, they have a nontrivial center).
To make them central, one has to restrict all functions to a
lattice, and so, to disconnect the leaf. 
After that we lose the effect of tunneling, but an interesting
picture of noncommutative geometry [17] is discovered.

{\bf Remark}.
The tunneling basic number $e^{-\pi^2/\hbar}$ 
appearing in formulas (7.11)--(7.15) admits a
geometrical interpretation via the membrane area. 
First of all, we remark that the remainder 
$O(e^{-\pi^2/\hbar})$ in (7.15) is bound up with 
the asymptotics of the probability operator $\cP$ (3.1). 
Formulas (3.15), (3.16) suggest an appropriate semiclassical
approximation of $\cP\sim\exp\{-\hbar\Delta/2\}$ 
as the heat operator with $\hbar$ playing the role
of ``time.''
But actually, in general, this is not correct
(look at Example 4.1).
Nevertheless, 
this approximation should work in the flat case.
Over the cylinder $\cX\approx\bR\times\bS$ with the usual flat
metric the heat kernel of $\exp(-\hbar\Delta/2)$ is given by 
\begin{align}
&\frac1{2\pi\hbar}\exp\Big\{-\frac{|z-z'|^2}{2\hbar}\Big\}
\theta\Big(\frac{i\pi(\oz-z-z'+\oz')}{\hbar},e^{-2\pi^2/\hbar}\Big)
\nonumber\\
&\qquad\qquad 
=
\frac1{2\pi\hbar}p_0(x,y)\big(1+O(e^{-2\pi^2/\hbar})\big),
\tag {7.19}
\end{align}
where $z=z(x)$ and $z'=z(y)$ are the complex coordinates 
of points $x,y\in\cX$. 
The exponents $2\pi^2j$ ($j=1,2,\dots$) appearing in this formula 
are the ``areas'' of nontrivial membranes
$\Sigma(x,x)\subset\cX^{\#}$  (see the notation before (2.10)).
Let us clarify: if the space $\cX$ is contractible,
then there are only trivial membranes $\Sigma(x,x)$ 
consisting of a single point, 
but if $\cX$ is not contractible, as in our case, 
then nontrivial membranes $\Sigma(x,x)$ exist.
For the first nontrivial membrane over the cylinder, we have
$\int_{\Sigma(x,x)}\omega_0^{\#}=2\pi^2$. 

We recall that the action $2\pi^2$ in (7.19) corresponds 
to the usual metric  over the cylinder. 
Our situation is more complicated 
since the quantum metric (7.13) differs from the usual one 
by the addition of order $O(e^{-\pi^2/\hbar})$.
This addition is generated by the reproducing kernel expansion
(7.12). 
The phase of the reproducing kernel is determined by triangle
membranes $\Sigma(x|y)$ which look like a ``half'' 
of the quadrangles $\Sigma(x,y)$.
Two sides of $\Sigma(x|y)$ belong to the $\pi_{\pm}$-fibers of
complex polarization, 
and the third side is the geodesic in $\cX$ between $x$ and $y$.
Let $\cX$ be not contractible.
Then even if $x=y$, we have nontrivial $\Sigma(x|x)$
corresponding to noncontractible closed geodesics
(or, in general, to noncontractible cycles on totally geodesic
Lagrangian submanifolds).
Over the cylinder, we see exactly 
$\int_{\Sigma(x|x)}\omega_0^{\#}=\pi^2$  
for the first nontrivial membrane $\Sigma(x|x)$. 
This is the geometrical explanation of the exponent in the
remainders $O(e^{-\pi^2/\hbar})$ in (7.12)--(7.14), 
and thus in the formula $p=p_0(1+O(e^{-\pi^2/\hbar}))$.
Since the Gauss probability function $p_0$ generates the
standard Wick product, we obtain (7.15) with the same remainder.

Following the sigma-model ideology, we could call $\Sigma(x,x)$
the mirror symmetric sigma-instanton, and 
$\Sigma(x|x)$ the asymmetric sigma-instanton.
Probably, these {\it sigma-instantons} generate additional ``complex 
stationary points''
in the ``path integral'' representation of the Kontsevich type
quantum products due to Cattaneo and Felder.

\section*{Acknowledgments}

This work was partially supported by RFBR (grant 99-01-01047).
The author also thanks Keio University, Dijon University, as
well as RIMS and IHES for their kind hospitality.  
Basically, the results of this work were presented 
at the Euroconference devoted to the memory of M.~Flato
(September 2000, Dijon). 
The author is grateful to 
A.~Alekseev, P.~Cartier, 
A.~Cattaneo, S.~Gutt, 
A.~Kirillov,
J.~Klauder, Sh.~Kobayashi, 
M.~Kontsevich, Y.~Kosmann-Schwarzbach, Y.~Maeda, G.~Marmo, 
T.~Miwa, H.~Omori, D.~Sternheimer, M.~Vergne, A.~Weinstein for 
important remarks. Also great thanks to E.~No\-vi\-ko\-va for
helpful discussions and collaboration.

\end{document}